\documentclass[11pt]{amsart}
\usepackage[margin=1.1in]{geometry}
\usepackage{amsmath,amssymb,amsthm,mathrsfs}
\usepackage{enumitem}
\usepackage{xcolor}
\usepackage[colorlinks=true,linkcolor=blue,citecolor=blue]{hyperref}
\usepackage[sc]{mathpazo}
\usepackage{microtype}
\usepackage{tikz}
\usetikzlibrary{arrows.meta,decorations.pathreplacing,patterns}
\usepackage[skins,breakable]{tcolorbox}
\usepackage{caption}
\usepackage{needspace}
\usepackage{etoolbox}
\BeforeBeginEnvironment{theorem}{\begin{samepage}}
\AfterEndEnvironment{theorem}{\end{samepage}}
\BeforeBeginEnvironment{lemma}{\begin{samepage}}
\AfterEndEnvironment{lemma}{\end{samepage}}
\BeforeBeginEnvironment{proposition}{\begin{samepage}}
\AfterEndEnvironment{proposition}{\end{samepage}}
\BeforeBeginEnvironment{corollary}{\begin{samepage}}
\AfterEndEnvironment{corollary}{\end{samepage}}
\BeforeBeginEnvironment{definition}{\begin{samepage}}
\AfterEndEnvironment{definition}{\end{samepage}}
\BeforeBeginEnvironment{example}{\Needspace{5\baselineskip}}
\definecolor{inkblue}{HTML}{1F4E9C}
\definecolor{calpha}{HTML}{3A56D4}
\definecolor{cbeta}{HTML}{C07A16}
\definecolor{cgamma}{HTML}{17876F}
\definecolor{ctau}{HTML}{7A4ED1}
\hypersetup{linkcolor=inkblue,citecolor=inkblue!70!black,urlcolor=inkblue}
\numberwithin{equation}{section}
\tikzset{
  gvert/.style={circle,draw=black!70,fill=black!4,inner sep=1.2pt,font=\small,minimum size=13pt},
  gvertsm/.style={circle,draw=black!60,fill=black!4,inner sep=0.9pt,font=\scriptsize,minimum size=10.5pt},
  gedge/.style={-{Stealth[length=2.1mm]},line width=.55pt,draw=black!60,shorten >=1pt,shorten <=1pt},
  gdim/.style={-{Stealth[length=1.6mm]},line width=.4pt,draw=black!22,shorten >=1pt,shorten <=1pt},
  ghl/.style={-{Stealth[length=2.2mm]},line width=1.05pt,shorten >=1pt,shorten <=1pt},
  ewt/.style={font=\tiny,text=black!55,inner sep=1.5pt}}
\newcommand\heapgrid[2]{%
  \foreach \i in {0,...,#1}{\draw[black!28,line width=.4pt,dash pattern=on .5pt off 1.5pt]
      ({\i*0.6},0) -- ({\i*0.6},{#2*0.46});}
  \foreach \i in {1,...,#1}{\node[font=\scriptsize,text=black!60] at ({(\i-0.5)*0.6},-0.26) {\i};}}
\newcommand\heappiece[5]{%
  \draw[rounded corners=2.2pt,fill=#4!20,draw=#4!60!black,line width=.7pt]
    ({(#1-1)*0.6+0.055},{(#3-1)*0.46+0.05}) rectangle ({#2*0.6-0.055},{#3*0.46-0.05});
  \node[font=\scriptsize,text=#4!35!black] at ({(#1-1+#2)*0.3},{(#3-0.5)*0.46}) {$#5$};}
\newcommand\heappieceK[5]{%
  \draw[rounded corners=2.2pt,fill=#4!14,draw=#4!75!black,line width=1pt]
    ({(#1-1)*0.6+0.055},{(#3-1)*0.46+0.05}) rectangle ({#2*0.6-0.055},{#3*0.46-0.05});
  \draw[rounded corners=2.2pt,pattern=north east lines,pattern color=#4!50]
    ({(#1-1)*0.6+0.055},{(#3-1)*0.46+0.05}) rectangle ({#2*0.6-0.055},{#3*0.46-0.05});
  \node[font=\scriptsize,text=#4!25!black,fill=white,fill opacity=.75,text opacity=1,inner sep=.6pt,rounded corners=1pt]
    at ({(#1-1+#2)*0.3},{(#3-0.5)*0.46}) {$#5$};}
\newcommand\heapword[3]{\node[font=\scriptsize,text=black!75] at ({#1*0.3},-0.72) {$#2$};\ignorespaces}
\newcommand\cslot[2]{\draw[black!40,line width=.45pt,dash pattern=on 1pt off 1.3pt]
   ({(#1-1)*0.42+0.045},{(#2-1)*0.42+0.045}) rectangle ({#1*0.42-0.045},{#2*0.42-0.045});}
\newcommand\ccellf[5]{\draw[fill=#3,draw=#4,line width=.6pt,rounded corners=1.3pt]
   ({(#1-1)*0.42+0.045},{(#2-1)*0.42+0.045}) rectangle ({#1*0.42-0.045},{#2*0.42-0.045});
   \node[font=\tiny] at ({(#1-0.5)*0.42},{(#2-0.5)*0.42}) {#5};}
\newcommand\ccell[4]{\ccellf{#1}{#2}{#3!22}{#3!60!black}{\textcolor{#3!30!black}{$#4$}}}

\newcommand\cjoin[5]{\draw[#5!50,line width=1.1pt] ({(#1-0.5)*0.42},{(#2-0.5)*0.42}) -- ({(#3-0.5)*0.42},{(#4-0.5)*0.42});}
\newcommand\ccols[1]{\foreach \i in {1,...,#1}{\node[font=\tiny,text=black!60] at ({(\i-0.5)*0.42},-0.19){\i};}}
\tcbset{
  mainthm/.style={enhanced jigsaw,breakable,colback=white,frame hidden,
    borderline west={2pt}{0pt}{inkblue!75},left=7pt,right=7pt,top=3pt,bottom=4pt,
    before skip=10pt plus 2pt,after skip=10pt plus 2pt},
  exbox/.style={enhanced jigsaw,breakable,colback=white,frame hidden,
    borderline west={1.4pt}{0pt}{black!30},left=7pt,right=7pt,top=2pt,bottom=4pt,
    before skip=9pt plus 2pt,after skip=9pt plus 2pt}}
\newcommand\heappieceKI[5]{%
  \draw[rounded corners=2.2pt,fill=#4!14,draw=#4!75!black,line width=1pt]
    ({(#1-1)*0.6+0.055},{(#3-1)*0.46+0.05}) rectangle ({#2*0.6-0.055},{#3*0.46-0.05});
  \draw[rounded corners=2.2pt,pattern=north west lines,pattern color=#4!50]
    ({(#1-1)*0.6+0.055},{(#3-1)*0.46+0.05}) rectangle ({#2*0.6-0.055},{#3*0.46-0.05});
  \node[font=\scriptsize,text=#4!25!black,fill=white,fill opacity=.75,text opacity=1,inner sep=.6pt,rounded corners=1pt]
    at ({(#1-1+#2)*0.3},{(#3-0.5)*0.46}) {$#5$};}
\newcommand\scell[5]{\draw[fill=#3,draw=#4,line width=.6pt,rounded corners=1.3pt]
   ({(#1-1)*0.5+0.05},{(#2-1)*0.42+0.04}) rectangle ({#1*0.5-0.05},{#2*0.42-0.04});
   \node[font=\scriptsize] at ({(#1-0.5)*0.5},{(#2-0.5)*0.42}) {#5};}
\newcommand\scols[2]{\foreach \i in {1,...,#1}{\node[font=\tiny,text=black!60] at ({(\i-0.5)*0.5},#2){\i};}}
\newcommand\lvl[3]{\node[font=\tiny,text=black!60,anchor=west] at ({#1*0.6+0.05},{(#2-0.5)*0.46}) {#3};}
\newcommand\membrow[5]{%
  \node[font=\tiny,anchor=east,text=black!60] at (-0.04,#3) {#4};
  \foreach \i in {1,...,#2}{\draw[black!40,line width=.4pt] ({(\i-0.5)*#1},#3) circle (0.045);}
  \foreach \i in {#5}{\fill[black!70] ({(\i-0.5)*#1},#3) circle (0.055);}}
\newtheorem{theorem}{Theorem}[section]
\newtheorem{lemma}[theorem]{Lemma}
\newtheorem{proposition}[theorem]{Proposition}
\newtheorem{corollary}[theorem]{Corollary}
\theoremstyle{definition}
\newtheorem{definition}[theorem]{Definition}
\newtheorem{example}[theorem]{Example}
\newtheorem{remark}[theorem]{Remark}
\newtheorem*{theoremA}{Leveled Heap Formula}
\newtheorem*{theoremB}{Positive Pair Rule}
\newtheorem*{theoremC}{Positivity Theorems}

\newcommand{\NN}{\mathbb{N}}
\newcommand{\ZZ}{\mathbb{Z}}
\newcommand{\QQ}{\mathbb{Q}}
\newcommand{\cA}{\mathcal{A}}
\newcommand{\cC}{\mathcal{C}}

\newcommand{\cN}{\mathcal{N}}
\newcommand{\cL}{\mathcal{L}}
\newcommand{\one}{\mathbf{1}}
\newcommand{\bH}{\mathbf{H}}
\newcommand{\supp}{\operatorname{supp}}
\newcommand{\cont}{\operatorname{cont}}
\newcommand{\Heap}{\operatorname{Heap}}
\newcommand{\SCC}{\operatorname{SCC}}
\newcommand{\Fun}{\operatorname{Fun}}
\newcommand{\Acyc}{\operatorname{Acyc}}
\newcommand{\Ch}{\operatorname{Ch}}
\newcommand{\wt}{\operatorname{wt}}

\title[Cycle popping and positivity]{Cycle popping and positivity in linear Laurent phenomenon algebras}
\author{Ross Glew}
\author{Pavlo Pylyavskyy}
\address{School of Mathematics, University of Minnesota, Minneapolis, MN 55455, USA}
\date{Draft, October 24, 2026}

\begin{document}

\begin{abstract}
We prove that every monomial in the cluster variables of a linear Laurent phenomenon algebra is a nonnegative combination of cluster monomials, as conjectured by Lam and the second author. The proof uses Wilson's cycle popping. We also give a new proof of Laurent positivity in linear Laurent phenomenon algebras, previously established by the first author. 
\end{abstract}

\maketitle

\section{Introduction}

\subsection{Cluster algebras and Laurent phenomenon algebras}

Cluster algebras were introduced by Fomin and Zelevinsky \cite{FZ1} to study total positivity and canonical bases in Lie theory. They have since appeared in many other areas of mathematics, including representation theory, Teichmüller theory, algebraic geometry, and integrable systems; see, for example, the survey \cite{Williams}. One of the central questions in the subject is positivity. Fomin and Zelevinsky conjectured that every cluster variable, when written as a Laurent polynomial in the variables of any cluster, has nonnegative coefficients. This conjecture was proved by Lee and Schiffler \cite{LeeSchiffler} for skew-symmetric cluster algebras and by Gross, Hacking, Keel, and Kontsevich \cite{GHKK} in general. One can also ask for a basis containing all cluster monomials in which multiplication has nonnegative structure constants. The theta bases of \cite{GHKK} provide such bases in broad generality.

Laurent phenomenon (LP) algebras were introduced by Lam and Pylyavskyy \cite{LP1} as a generalization of cluster algebras in which exchange polynomials are no longer required to be binomials. The Laurent phenomenon continues to hold in this setting. Positivity, however, does not hold for arbitrary LP algebras, so it becomes a question about particular families. In \cite{LP2}, Lam and Pylyavskyy studied a family of linear LP algebras associated with directed graphs. Their clusters have a particularly nice combinatorial description in terms of nested collections of strongly connected subsets, generalizing the associahedron and other graph associahedra. These algebras were subsequently studied in \cite{BCKZ1,BCKZ2,DTW}.

It was conjectured in \cite{LP2} that these linear LP algebras satisfy two positivity properties analogous to those of cluster algebras. First, every cluster variable is a Laurent polynomial with nonnegative coefficients in the variables of any cluster. Second, every monomial in cluster variables is a nonnegative linear combination of cluster monomials. For the second statement, the proposed positive basis is already given to us: the cluster monomials themselves form a basis \cite{DTW}. No additional basis elements are needed. The first conjecture was recently proved by Glew \cite{Glew}. In this paper we prove the second conjecture and give a new proof of the first. The two proofs use the same combinatorial identities.

\subsection{Linear Laurent phenomenon algebras and positivity conjectures}

Let $\Gamma$ be a directed graph on $[n]$. We recall the
linear LP algebra $\cA_\Gamma$ associated with it. Its initial cluster is $X_1,\ldots,X_n$, with coefficient variables
$A_1,\ldots,A_n$. The singleton cluster variables satisfy
\[
X_iY_i=A_i+\sum_{i\to j}X_j.
\]
The other cluster variables are indexed by nonempty strongly
connected subsets $I\subseteq[n]$. They are the principal minors
\[
Y_I=\det(D_Y-C)_{I,I},
\qquad D_Y=\operatorname{diag}(Y_1,\ldots,Y_n),
\]
where $C$ is the adjacency matrix of $\Gamma$ \cite{LP2}.
We use the same definition for arbitrary subsets $I$, with
$Y_\emptyset=1$. The minor $Y_I$ factors over the strongly
connected components of $\Gamma[I]$.

The clusters are described by nested collections of strongly
connected subsets, together with the appropriate initial
variables $X_i$. The cluster monomials form a basis of
$\cA_\Gamma$ \cite{DTW}. The two positivity conjectures concern different expansions:
Laurent expansions in a cluster, and multiplication in the
cluster-monomial basis. Neither follows from the Laurent phenomenon.

We will keep an independent weight $c_{ij}$ on each edge $i\to j$. This gives weighted principal minors and weighted versions of the $X_i$, and lets us see which cycles and paths contribute to each coefficient. All our combinatorial identities hold with these weights. At the end we set $c_{ij}=1$ to recover the original linear LP algebra; no mutation theory for the weighted system is needed.

\subsection{Statements of the main results}

The main idea is to write the principal minors in coordinates where they are positive sums. At each vertex $i\in S$, choose either a root of weight $A_i/X_i$ or an arrow $i\to j$ of weight $c_{ij}X_j/X_i$. Arrows may leave $S$, but they must not form a directed cycle inside $S$. Such a collection of choices is called an \emph{acyclic function}. The minor $Y_S$ is the sum of their weights (Proposition~\ref{prop:forest}).

A monomial $Y^T=\prod_iY_i^{T_i}$ therefore describes
stacks of arrows, with $T_i$ entries at vertex $i$.
We use Wilson's cycle-popping procedure \cite{Wilson,PW}:
pop cycles formed by the exposed arrows, remove the
remaining acyclic function as one layer, and repeat.
Each popped cycle contributes only its edge weight, since the $X$-variables telescope around the cycle. Disjoint exposed cycles can be popped in either order. We therefore record a history up to swaps of successive disjoint cycles, rather than as a single ordered sequence. This gives a heap in the sense of Cartier--Foata and Viennot \cite{CF,Viennot}.

To state our first result, for $U\in\NN^n$ let
\[
U^{(r)}=\{i:U_i\ge r\},
\qquad
B_U=\prod_{r\ge1}Y_{U^{(r)}}.
\]
The sets $U^{(1)}\supseteq U^{(2)}\supseteq\cdots$
are the layers of $U$. The elements $B_U$ are
the weighted versions of the cluster monomials
involving only $Y$-variables.

\begin{theoremA}[Theorem~\ref{thm:A}]
For $T\in\NN^n$,
\[
Y^T=\sum_U\left(\sum_{\bH}c^{\bH}\right)B_U.
\]
Here $\bH=(H_1,H_2,\ldots)$ runs over sequences
of heaps such that
\[
H_r\in\Heap(U^{(r-1)}),
\qquad
\sum_r\cont(H_r)=T-U,
\]
where $U^{(0)}=[n]$.
\end{theoremA}

The index $r$ counts popping rounds: $H_r$ is popped after $r-1$ acyclic layers have been removed. It is not the height of a piece within a heap. The same cycles may occur in different rounds, and these possibilities account for the multiplicities in Theorem~\ref{thm:A}.

Our second result gives a positive rule for multiplying
two principal minors, analogous to a
Littlewood--Richardson rule.

\begin{theoremB}[Theorem~\ref{thm:B}]
For $I,J\subseteq[n]$ and $\kappa=\one_I+\one_J$,
\[
Y_IY_J=\sum_H c^H B_{\kappa-\cont(H)}.
\]
The sum is over heaps with $\cont(H)\le\kappa$
such that no maximal piece is supported in $I$
and no minimal piece is supported in $J$.
\end{theoremB}

To prove the Positive Pair Rule (Theorem~\ref{thm:B}), we need only the case of the Leveled Heap Formula (Theorem~\ref{thm:A}) in which each vertex occurs at most twice. We expand $Y^\kappa$ in two ways: separately in its two factors, and all at once. A finite bijection strips pieces from the two boundaries of a heap and identifies the contributions for smaller pairs of sets. Induction then gives the formula.

The pair rule has exactly the form needed for straightening. Its empty-heap term replaces $I,J$ by the more nested pair $I\cup J,I\cap J$, while every nonempty heap lowers total vertex multiplicity. Repeated application therefore gives positive multiplication in the basis $\{B_U\}$. An exit-path identity handles the products involving $X$-variables. Together these give our main positivity results.

\begin{theoremC}[Corollaries~\ref{cor:LP-positivity} and~\ref{cor:laurent-LP}]
For every directed graph $\Gamma$:
\begin{enumerate}[leftmargin=2.2em,label=\textup{(\roman*)}]
\item Every monomial in the cluster variables of
$\cA_\Gamma$ is a linear combination of cluster
monomials with coefficients in
$\NN[A_1,\ldots,A_n]$.
\item Every cluster variable has a Laurent expansion
with coefficients in $\NN[A_1,\ldots,A_n]$
in every cluster.
\end{enumerate}
\end{theoremC}

Part~(i) proves \cite[Conjecture~7.3(ii)]{LP2}.
Dantas e Moura, Telekicherla Kandalam, and Woodruff
previously proved it for bidirected trees, as well
as positive product formulas for disjoint sets
and sets meeting in one vertex \cite{DTW}.

For Laurent positivity, Banaian, Chepuri, Kelley,
and Zhang gave formulas for rooted tree clusters
\cite{BCKZ1} and further tree cases \cite{BCKZ2}.
The general result was recently proved by the
first author \cite{Glew}. Our proof follows the
same broad strategy, but uses Theorem~\ref{thm:B}
in place of the iterated path identity, and organizes
the induction directly by the nested structure of the
cluster. Iterating the recursion gives explicit positive
Laurent expressions.

The paper is organized as follows.
Section~\ref{sec:arrows} introduces acyclic
functions and the basis $B_U$.
Section~\ref{sec:popping} proves the {Leveled Heap Formula}
by cycle popping, and Section~\ref{sec:pair} derives the
{Positive Pair Rule} by stripping heap boundaries.
Positive straightening and Laurent positivity are developed
in Sections~\ref{sec:straightening} and~\ref{sec:laurent}.

\subsection{A small example}

Consider the directed graph
\[
1\leftrightarrow2\leftrightarrow3\leftrightarrow4,
\qquad 3\longrightarrow1,
\]
which will serve as our running example.
Write $\alpha=(12)$, $\beta=(23)$ and $\tau=(123)$
for three of its directed cycles, with weights
$f_{12}=c_{12}c_{21}$, $f_{23}=c_{23}c_{32}$,
and $t=c_{12}c_{23}c_{31}$. Then Theorem~\ref{thm:B} gives
\begin{equation}\label{eq:intro}
\boxed{
Y_{12}Y_{23}=Y_{123}Y_2+tY_2+f_{12}f_{23}.
}
\end{equation}

The three terms correspond to the empty heap,
the triangle $\tau$, and the heap with $\alpha$
below $\beta$; see Figure~\ref{fig:intro}.
The product on the left has one occurrence each of vertices $1$ and $3$, and two of vertex $2$. Each cycle uses one occurrence of every vertex it meets. Reading the unused occurrences in layers gives the minor factors on the right: the empty heap leaves $123$ and $2$, the triangle leaves just $2$, and the two overlapping cycles use everything.

\par\medskip
\begin{center}
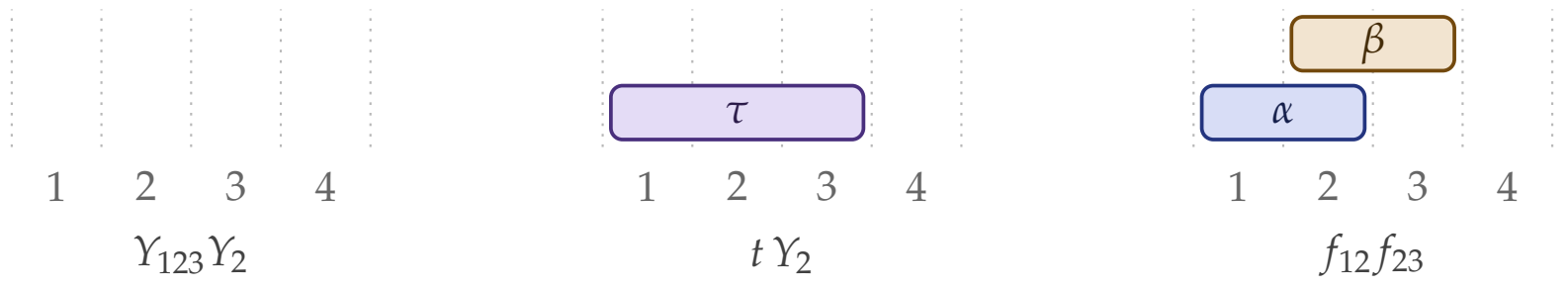
\begin{minipage}{\linewidth}\centering
\begin{tikzpicture}[baseline=0pt]
  \heapgrid{4}{2}
  \heapword{4}{Y_{123}Y_2}{}
\end{tikzpicture}\qquad\qquad
\begin{tikzpicture}[baseline=0pt]
  \heapgrid{4}{2}
  \heappiece{1}{3}{1}{ctau}{\tau}
  \heapword{4}{t\,Y_2}{}
\end{tikzpicture}\qquad\qquad
\begin{tikzpicture}[baseline=0pt]
  \heapgrid{4}{2}
  \heappiece{1}{2}{1}{calpha}{\alpha}
  \heappiece{2}{3}{2}{cbeta}{\beta}
  \heapword{4}{f_{12}f_{23}}{}
\end{tikzpicture}
\captionof{figure}{The three heaps indexing the terms
of \eqref{eq:intro}. Columns are vertices; each bar is
a cycle, resting on the bars it meets. No top piece
lies inside $\{1,2\}$ and no bottom piece lies
inside $\{2,3\}$.}\label{fig:intro}
\end{minipage}\end{center}
\medskip

\subsection*{Conventions}
We use $\NN=\{0,1,2,\ldots\}$.
For $T\in\NN^n$, write $|T|=\sum_iT_i$ and
$\supp T=\{i:T_i>0\}$.
Inequalities between such vectors are
coordinatewise. The vector $\one_S$ is
the indicator of $S$, and
$Y^T=\prod_iY_i^{T_i}$.
Empty products have value $1$.
We abbreviate vertex sets in subscripts:
for example, $Y_{123}=Y_{\{1,2,3\}}$.

\subsection{AI statement} ChatGPT (OpenAI) and Claude (Anthropic) assisted with mathematical exploration, verification, and preparation of the manuscript. The cycle-popping approach and the strategy for Laurent positivity originated with the authors, while the boundary-stripping proof of the pair rule was developed with substantial AI assistance. The authors take full responsibility for the content.

\subsection{Overlapping work}
While finishing this paper, we learned of a recent preprint
by Tang \cite{Tang}. His Theorem~8.2 also establishes
positivity of cluster-monomial multiplication, as part of a
stronger atomicity result, using the Laurent positivity theorem
of the first author \cite{Glew} and a monomial-separation argument.
Our approach is different and more constructive: it gives explicit
positive formulas for products of variables
(Theorem~\ref{thm:B}), as well as positive Laurent expansions
obtained through cycle popping.

\section{Acyclic functions and principal minors}\label{sec:arrows}

\subsection{Weighted principal minors}
Let $\Gamma$ be a finite directed graph on $[n]$, without loops or multiple edges. Both $i\to j$ and $j\to i$ may occur. Give each edge an indeterminate weight $c_{ij}$, and write $\ZZ[c]$ for the polynomial ring in these weights. Let $Y_1,\ldots,Y_n$ be indeterminates and set
\[
Y_S=\det(D_Y-C)_{S,S},\qquad Y_\emptyset=1,
\]
where $C_{ij}=c_{ij}$ for an edge $i\to j$, and $C_{ij}=0$ otherwise. At $c_{ij}=1$, the minors indexed by nonempty strongly connected sets are the $Y$-cluster variables of $\cA_\Gamma$ \cite[Proposition~4.5]{LP2}.

We will also use $Y_S$ when $S$ is not strongly connected. Ordering the strongly connected components of $\Gamma[S]$ topologically makes the matrix block triangular, so
\begin{equation}\label{eq:scc}
Y_S=\prod_{K\in\SCC(\Gamma[S])}Y_K.
\end{equation}
This lets us use arbitrary vertex sets and separate their strongly connected components only when needed.

A \emph{cycle} is a directed simple cycle of length at least two, considered up to cyclic rotation. Opposite orientations are different cycles when the length is at least three. For a cycle $\gamma$, write $V(\gamma)$ for its vertex set and $c^\gamma$ for the product of its edge weights. Let $\cC(S)$ be the cycles supported in $S$. A \emph{clique} is a collection $K$ of pairwise vertex-disjoint cycles, including the empty collection. Put
\[
V(K)=\bigcup_{\gamma\in K}V(\gamma),\qquad
c^K=\prod_{\gamma\in K}c^\gamma.
\]

\begin{example}[The running graph]\label{ex:running}
We use the graph
\[
1\leftrightarrow2\leftrightarrow3\leftrightarrow4,
\qquad 3\longrightarrow1
\]
throughout. It has four cycles,
\[
\alpha=(12),\quad\beta=(23),\quad\gamma=(34),\quad\tau=(123),
\]
with weights
\[
f_{12}=c_{12}c_{21},\quad f_{23}=c_{23}c_{32},\quad
f_{34}=c_{34}c_{43},\quad t=c_{12}c_{23}c_{31}.
\]
Only $\alpha$ and $\gamma$ are disjoint.

\par\medskip
\begin{center}\begin{minipage}{\linewidth}\centering
\begin{tikzpicture}
  \node[gvert] (v1) at (0,0) {$1$};
  \node[gvert] (v2) at (1.65,0) {$2$};
  \node[gvert] (v3) at (3.3,0) {$3$};
  \node[gvert] (v4) at (4.95,0) {$4$};
  \draw[gedge] (v1) to[bend left=22] node[ewt,above]{$c_{12}$} (v2);
  \draw[gedge] (v2) to[bend left=22] node[ewt,below]{$c_{21}$} (v1);
  \draw[gedge] (v2) to[bend left=22] node[ewt,above]{$c_{23}$} (v3);
  \draw[gedge] (v3) to[bend left=22] node[ewt,below]{$c_{32}$} (v2);
  \draw[gedge] (v3) to[bend left=22] node[ewt,above]{$c_{34}$} (v4);
  \draw[gedge] (v4) to[bend left=22] node[ewt,below]{$c_{43}$} (v3);
  \draw[gedge] (v3) to[bend left=48] node[ewt,below]{$c_{31}$} (v1);
\end{tikzpicture}\\[11pt]
\begin{tabular}{@{}c@{\hspace{14pt}}c@{\hspace{14pt}}c@{\hspace{14pt}}c@{}}
\begin{tikzpicture}[scale=0.62]
  \node[gvertsm] (w1) at (0.00,0) {$1$};
  \node[gvertsm] (w2) at (1.42,0) {$2$};
  \node[gvertsm] (w3) at (2.84,0) {$3$};
  \node[gvertsm] (w4) at (4.26,0) {$4$};
  \draw[gdim] (w1) to[bend left=22] (w2);
  \draw[gdim] (w2) to[bend left=22] (w1);
  \draw[gdim] (w2) to[bend left=22] (w3);
  \draw[gdim] (w3) to[bend left=22] (w2);
  \draw[gdim] (w3) to[bend left=22] (w4);
  \draw[gdim] (w4) to[bend left=22] (w3);
  \draw[gdim] (w3) to[bend left=48] (w1);
  \draw[ghl,draw=calpha] (w1) to[bend left=22] (w2);
  \draw[ghl,draw=calpha] (w2) to[bend left=22] (w1);
  \node[font=\scriptsize,text=calpha!30!black] at (2.13,-0.78) {$\alpha=(12)$};
\end{tikzpicture} &
\begin{tikzpicture}[scale=0.62]
  \node[gvertsm] (w1) at (0.00,0) {$1$};
  \node[gvertsm] (w2) at (1.42,0) {$2$};
  \node[gvertsm] (w3) at (2.84,0) {$3$};
  \node[gvertsm] (w4) at (4.26,0) {$4$};
  \draw[gdim] (w1) to[bend left=22] (w2);
  \draw[gdim] (w2) to[bend left=22] (w1);
  \draw[gdim] (w2) to[bend left=22] (w3);
  \draw[gdim] (w3) to[bend left=22] (w2);
  \draw[gdim] (w3) to[bend left=22] (w4);
  \draw[gdim] (w4) to[bend left=22] (w3);
  \draw[gdim] (w3) to[bend left=48] (w1);
  \draw[ghl,draw=cbeta] (w2) to[bend left=22] (w3);
  \draw[ghl,draw=cbeta] (w3) to[bend left=22] (w2);
  \node[font=\scriptsize,text=cbeta!30!black] at (2.13,-0.78) {$\beta=(23)$};
\end{tikzpicture} &
\begin{tikzpicture}[scale=0.62]
  \node[gvertsm] (w1) at (0.00,0) {$1$};
  \node[gvertsm] (w2) at (1.42,0) {$2$};
  \node[gvertsm] (w3) at (2.84,0) {$3$};
  \node[gvertsm] (w4) at (4.26,0) {$4$};
  \draw[gdim] (w1) to[bend left=22] (w2);
  \draw[gdim] (w2) to[bend left=22] (w1);
  \draw[gdim] (w2) to[bend left=22] (w3);
  \draw[gdim] (w3) to[bend left=22] (w2);
  \draw[gdim] (w3) to[bend left=22] (w4);
  \draw[gdim] (w4) to[bend left=22] (w3);
  \draw[gdim] (w3) to[bend left=48] (w1);
  \draw[ghl,draw=cgamma] (w3) to[bend left=22] (w4);
  \draw[ghl,draw=cgamma] (w4) to[bend left=22] (w3);
  \node[font=\scriptsize,text=cgamma!30!black] at (2.13,-0.78) {$\gamma=(34)$};
\end{tikzpicture} &
\begin{tikzpicture}[scale=0.62]
  \node[gvertsm] (w1) at (0.00,0) {$1$};
  \node[gvertsm] (w2) at (1.42,0) {$2$};
  \node[gvertsm] (w3) at (2.84,0) {$3$};
  \node[gvertsm] (w4) at (4.26,0) {$4$};
  \draw[gdim] (w1) to[bend left=22] (w2);
  \draw[gdim] (w2) to[bend left=22] (w1);
  \draw[gdim] (w2) to[bend left=22] (w3);
  \draw[gdim] (w3) to[bend left=22] (w2);
  \draw[gdim] (w3) to[bend left=22] (w4);
  \draw[gdim] (w4) to[bend left=22] (w3);
  \draw[gdim] (w3) to[bend left=48] (w1);
  \draw[ghl,draw=ctau] (w1) to[bend left=22] (w2);
  \draw[ghl,draw=ctau] (w2) to[bend left=22] (w3);
  \draw[ghl,draw=ctau] (w3) to[bend left=48] (w1);
  \node[font=\scriptsize,text=ctau!30!black] at (2.13,-0.78) {$\tau=(123)$};
\end{tikzpicture}
\end{tabular}

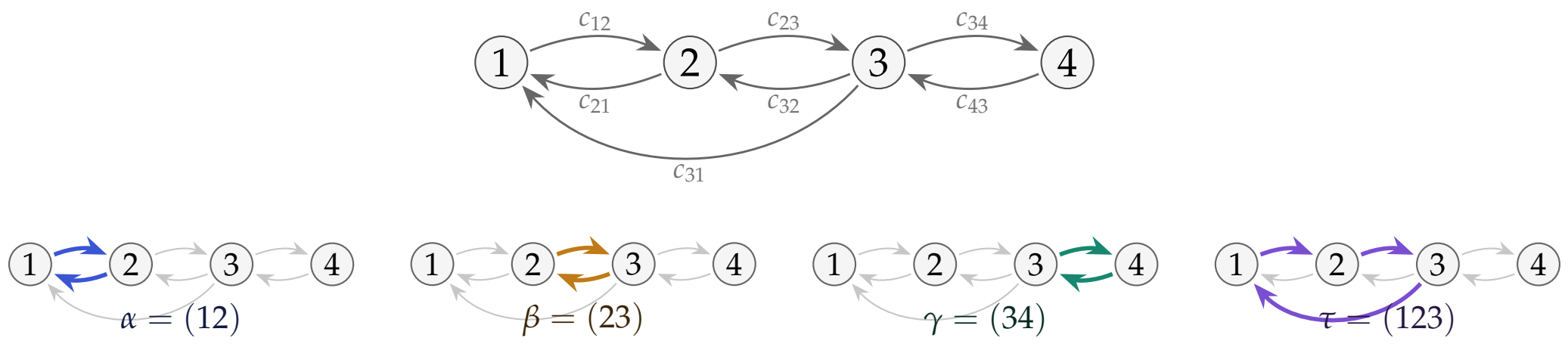
\captionof{figure}{The running example: the graph $\Gamma$ with its edge weights (top), and its four directed cycles (bottom).}\label{fig:running-graph}
\end{minipage}\end{center}
\medskip

The first few minors are
\[
\begin{aligned}
Y_{12}&=Y_1Y_2-f_{12},\\
Y_{123}&=Y_1Y_2Y_3-f_{12}Y_3-f_{23}Y_1-t,\\
Y_{1234}&=Y_1Y_2Y_3Y_4-f_{12}Y_3Y_4-f_{23}Y_1Y_4
-f_{34}Y_1Y_2-tY_4+f_{12}f_{34}.
\end{aligned}
\]
Notice that the positive term $f_{12}f_{34}$ comes from the only pair of disjoint cycles, $\alpha,\gamma$. The induced graph on $124$ has components $12$ and $4$, so $Y_{124}=Y_{12}Y_4$.
\end{example}

The determinant expansion gives, in general,
\begin{equation}\label{eq:perm}
Y_S=\sum_{K\subseteq\cC(S)\text{ a clique}}
(-1)^{|K|}c^K Y^{\one_{S\setminus V(K)}}.
\end{equation}
Indeed, the nontrivial cycles of a permutation form $K$. A cycle of length $\ell$ has permutation sign $(-1)^{\ell-1}$ and contributes $\ell$ negative matrix entries, so its total contribution is $-c^\gamma$. 

\subsection{Arrows and acyclic functions}
The determinant formula identifies the variables, but its signs do not yet explain positivity. For this we return to the initial-cluster coordinates.

Introduce independent variables $X_1,\ldots,X_n,A_1,\ldots,A_n$ and work in
\[
F=\QQ(c)(X_1,\ldots,X_n,A_1,\ldots,A_n).
\]
Make the substitution
\begin{equation}\label{eq:Yi}
Y_i=\frac{A_i+\sum_{i\to j}c_{ij}X_j}{X_i}.
\end{equation}
This is simply a change of variables: keeping the $X_i$ fixed,
we can recover the $A_i$ from the $Y_i$ by
\[
A_i=X_iY_i-\sum_j c_{ij}X_j.
\]
Thus the $Y_i$ remain algebraically independent, and any
identity proved using the arrow model holds for our original
principal minors. At $c_{ij}=1$, \eqref{eq:Yi} is the exchange
relation of the linear LP algebra.

At a vertex $v$, choose either a root $*$ or an allowed arrow $v\to j$. Their weights are
\[
\wt_v(*)=\frac{A_v}{X_v},\qquad
\wt_v(j)=\frac{c_{vj}X_j}{X_v}.
\]
Write $\Ch(v)$ for this set of choices. A \emph{function on $S$} makes one choice at each $v\in S$; thus
\[
\Fun(S)=\prod_{v\in S}\Ch(v),\qquad
\wt(\varphi)=\prod_{v\in S}\wt_v(\varphi(v)).
\]
Its directed cycles inside $S$ are necessarily vertex-disjoint. We denote their clique by $K(\varphi)$ and call $\varphi$ \emph{acyclic} when $K(\varphi)=\emptyset$. Let $\Acyc(S)$ be the set of acyclic functions on $S$.

An arrow leaving $S$ is allowed. We stop following a path when it leaves $S$ or reaches a root; only cycles entirely inside $S$ are forbidden. Around such a cycle, the $X$-variables cancel:
\begin{equation}\label{eq:telescope}
\prod_{v\to j\text{ in }\gamma}\frac{c_{vj}X_j}{X_v}=c^\gamma.
\end{equation}

\begin{proposition}[Acyclic-function formula]\label{prop:forest}
For every $S\subseteq[n]$,
\[
Y_S=\sum_{\varphi\in\Acyc(S)}\wt(\varphi).
\]
\end{proposition}

The unweighted formula was proved in the opposite direction
in \cite[Proposition~4.5]{LP2}: there the $Y_S$ were defined
by acyclic functions and identified with principal minors.
Here we start from the minors and allow arbitrary edge weights.
This is also a form of the directed matrix-tree theorem
\cite{Chaiken}.

\begin{proof}
Expand each singleton factor in \eqref{eq:perm} using \eqref{eq:Yi}. A term consists of a clique $K$ and a function on $S\setminus V(K)$. Adding the arrows of $K$ gives a function $\varphi$ on $S$, and the possible choices of $K$ are exactly the subcollections of $K(\varphi)$. By \eqref{eq:telescope}, their total contribution is
\[
\wt(\varphi)\sum_{K\subseteq K(\varphi)}(-1)^{|K|}.
\]
This is $\wt(\varphi)$ when $\varphi$ is acyclic and zero otherwise.
\end{proof}

A product of singleton variables counts all functions. Separating them according to their cycles gives a positive identity in the opposite direction.

\begin{corollary}[Removing the cycles of a function]\label{cor:cycdec}
Removing the cycles gives a weight-preserving bijection
\[
\Fun(S)\longleftrightarrow
\bigsqcup_{K\subseteq\cC(S)\text{ a clique}}
\{K\}\times\Acyc(S\setminus V(K)),
\]
where the weight on the right is $c^K\wt(\psi)$. Consequently,
\begin{equation}\label{eq:cycdec}
Y^{\one_S}=\sum_{K\subseteq\cC(S)\text{ a clique}}c^K Y_{S\setminus V(K)}.
\end{equation}
\end{corollary}

\begin{proof}
The restriction outside all cycle vertices is acyclic. Conversely, adjoining a clique $K$ to an acyclic function on its complement creates no extra cycle: once an arrow path enters $V(K)$, it stays on one of the cycles of $K$. The weights agree by \eqref{eq:telescope}.
\end{proof}

\begin{example}\label{ex:forest}
There are six functions on $12$ in the running graph: vertex $1$ chooses $*$ or $2$, and vertex $2$ chooses $*$, $1$, or $3$. Only the choice $(1\mapsto2,2\mapsto1)$ is cyclic. Its weight is $f_{12}$, so
\[
Y_1Y_2=Y_{12}+f_{12}.
\]
Figure~\ref{fig:forest} illustrates roots, exits, and a cycle on the slightly larger set $123$.

\par\medskip
\begin{center}\begin{minipage}{\linewidth}\centering
\newcommand\fgraph{%
  \node[gvertsm] (w1) at (0.00,0) {$1$};
  \node[gvertsm] (w2) at (1.32,0) {$2$};
  \node[gvertsm] (w3) at (2.64,0) {$3$};
  \node[gvertsm] (w4) at (3.96,0) {$4$};
  \draw[gdim] (w1) to[bend left=20] (w2); \draw[gdim] (w2) to[bend left=20] (w1);
  \draw[gdim] (w2) to[bend left=20] (w3); \draw[gdim] (w3) to[bend left=20] (w2);
  \draw[gdim] (w3) to[bend left=20] (w4); \draw[gdim] (w4) to[bend left=20] (w3);
  \draw[gdim] (w3) to[bend left=40] (w1);}%
\begin{tabular}{@{}c@{\hspace{14pt}}c@{\hspace{14pt}}c@{}}
\begin{tikzpicture}[scale=0.82,baseline=(w1.base)]
  \fill[inkblue!6,rounded corners=8pt] (-0.42,-0.95) rectangle (3.06,0.6);
  \fgraph
  \draw[ghl,draw=inkblue!85] (w1) to[bend left=20] (w2);
  \draw[ghl,draw=inkblue!85] (w3) to[bend left=20] (w4);
  \draw[inkblue!85,line width=.95pt] (w2.center) circle (0.25);
  \node[font=\small] at (1.95,-1.55) {$\dfrac{c_{12}c_{34}\,A_2X_4}{X_1X_3}$};
\end{tikzpicture} &
\begin{tikzpicture}[scale=0.82,baseline=(w1.base)]
  \fill[inkblue!6,rounded corners=8pt] (-0.42,-0.95) rectangle (3.06,0.6);
  \fgraph
  \draw[ghl,draw=inkblue!85] (w1) to[bend left=20] (w2);
  \draw[ghl,draw=inkblue!85] (w3) to[bend left=20] (w2);
  \draw[inkblue!85,line width=.95pt] (w2.center) circle (0.25);
  \node[font=\small] at (1.95,-1.55) {$\dfrac{c_{12}c_{32}\,A_2X_2}{X_1X_3}$};
\end{tikzpicture} &
\begin{tikzpicture}[scale=0.82,baseline=(w1.base)]
  \fill[red!5,rounded corners=8pt] (-0.42,-0.95) rectangle (3.06,0.6);
  \fgraph
  \draw[ghl,draw=ctau] (w1) to[bend left=20] (w2);
  \draw[ghl,draw=ctau] (w2) to[bend left=20] (w3);
  \draw[ghl,draw=ctau] (w3) to[bend left=40] (w1);
  \node[font=\small] at (1.95,-1.55) {cycle $\tau$:\ \ $t$};
\end{tikzpicture}
\end{tabular}

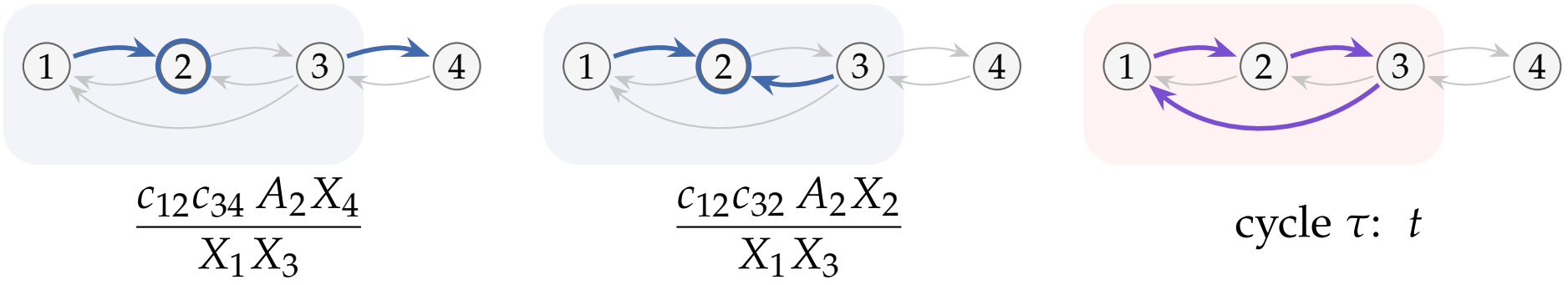
\captionof{figure}{Three functions on $S=\{1,2,3\}$ (shaded), with their weights; a circled vertex is a root. Left: a root at $2$ and an arrow $3\to4$ leaving $S$; this function is acyclic. Middle: an acyclic function whose arrows all lead to the root $2$. Right: the arrows close the cycle $\tau$ inside $S$, so this function is not acyclic; its weight telescopes to $c^\tau=t$.}\label{fig:forest}
\end{minipage}\end{center}
\medskip

\end{example}

\subsection{Layers and the basis}
We will write the positive expansions in a basis of nested products. A convenient way to index these products is by their vertex multiplicities. Read $U\in\NN^n$ as a diagram with $U_i$ boxes in column $i$. Its rows are
\[
U^{(r)}=\{i:U_i\ge r\},\qquad r\ge1.
\]
They form a decreasing sequence, and $U=\sum_{r\ge1}\one_{U^{(r)}}$. We also set $U^{(0)}=[n]$. Define
\[
B_U=\prod_{r\ge1}Y_{U^{(r)}}.
\]
By \eqref{eq:scc}, each layer factors into strongly connected components. These components form the unique nested multicollection with total multiplicity $U$ \cite[Lemma~2.1]{DTW}. Thus, at $c_{ij}=1$, the $B_U$ are exactly the cluster monomials involving only $Y$-variables.

The determinant has leading term $Y^{\one_S}$. Hence
\[
B_U=Y^U+\text{terms of smaller total degree in the }Y_i,
\]
and the $B_U$ form a $\ZZ[c]$-basis of $\ZZ[c][Y_1,\ldots,Y_n]$.

\begin{example}\label{ex:layers}
For $U=(3,2,1,2)$, the layers are $1234$, $124$, and $1$. Therefore
\[
B_U=Y_{1234}Y_{124}Y_1
=Y_{1234}Y_{12}Y_4Y_1.
\]
The factor $Y_4$ comes from the second layer, not from a new row of the diagram.
\end{example}

\section{Cycle popping and the leveled heap formula}\label{sec:popping}

\subsection{Heaps of cycles}
Heaps will let us remember the order of overlapping cycles without choosing an order for independent removals.

A \emph{heap} is a word in directed cycles, modulo the relations
\[
\gamma\delta=\delta\gamma
\qquad\text{when }V(\gamma)\cap V(\delta)=\emptyset.
\]
The pieces of a heap are \emph{occurrences} of cycles: two copies of the same cycle are two different pieces. To read the order from a representative word, put an earlier piece below a later one whenever their supports meet, then take the transitive closure. Swapping adjacent disjoint pieces does not change this partial order. Every pair of pieces sharing a vertex is comparable, and every covering relation joins pieces whose supports meet.

The minimal pieces are pairwise disjoint, as are the maximal pieces. Multiplication is stacking: in $H\cdot H'$, the pieces of $H'$ are placed above those of $H$. Write
\[
c^H=\prod_{\gamma\text{ a piece of }H}c^\gamma,
\qquad
\cont(H)_v=\#\{\text{pieces of }H\text{ through }v\}.
\]
Weight is multiplicative and content is additive. A heap lives in $S$ when every piece is supported in $S$; the set of such heaps is $\Heap(S)$.

\begin{example}\label{ex:heaps}
For the running graph, the only commuting pair is $\alpha,\gamma$. Thus $\alpha\beta$ and $\beta\alpha$ are different heaps, whereas $\alpha\gamma=\gamma\alpha$. Disjoint pieces need not be incomparable: in $\alpha\beta\gamma$, the piece $\beta$ forces $\alpha$ below $\gamma$, even though those two cycles are disjoint.

\par\medskip
\begin{center}
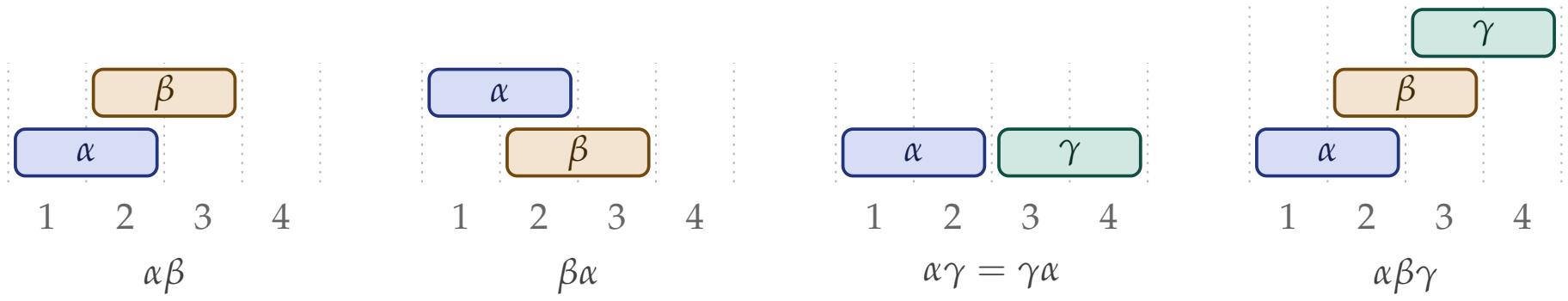
\begin{minipage}{\linewidth}\centering
\begin{tikzpicture}[baseline=0pt]
  \heapgrid{4}{2}
  \heappiece{1}{2}{1}{calpha}{\alpha}
  \heappiece{2}{3}{2}{cbeta}{\beta}
  \heapword{4}{\alpha\beta}{}
\end{tikzpicture}\qquad
\begin{tikzpicture}[baseline=0pt]
  \heapgrid{4}{2}
  \heappiece{2}{3}{1}{cbeta}{\beta}
  \heappiece{1}{2}{2}{calpha}{\alpha}
  \heapword{4}{\beta\alpha}{}
\end{tikzpicture}\qquad
\begin{tikzpicture}[baseline=0pt]
  \heapgrid{4}{2}
  \heappiece{1}{2}{1}{calpha}{\alpha}
  \heappiece{3}{4}{1}{cgamma}{\gamma}
  \heapword{4}{\alpha\gamma=\gamma\alpha}{}
\end{tikzpicture}\qquad
\begin{tikzpicture}[baseline=0pt]
  \heapgrid{4}{3}
  \heappiece{1}{2}{1}{calpha}{\alpha}
  \heappiece{2}{3}{2}{cbeta}{\beta}
  \heappiece{3}{4}{3}{cgamma}{\gamma}
  \heapword{4}{\alpha\beta\gamma}{}
\end{tikzpicture}
\captionof{figure}{Heaps of cycles of the running example. In every cycle of this graph the vertices form an interval, so each piece is drawn as one bar.}\label{fig:heaps}
\end{minipage}\end{center}
\medskip

\end{example}

\subsection{Stacks and cycle popping}
A \emph{stack configuration of height $T$} is a word
$\sigma_v=(\sigma_v(1),\ldots,\sigma_v(T_v))$ in $\Ch(v)$ at each vertex $v$. The first entry is the top. The weight of the configuration is the product of all entry weights, so
\begin{equation}\label{eq:stacks}
Y^T=\sum_{\sigma\text{ of height }T}\wt(\sigma).
\end{equation}
The top entries form the \emph{exposed function} on the vertices with nonempty stacks. If some of these entries form a cycle $\gamma$, we may \emph{pop} it by removing the top entry at every vertex of $\gamma$. The removed entries have total weight $c^\gamma$, so the original weight is $c^\gamma$ times the weight of what remains.

Continue until the exposed function is acyclic. The process is finite, since each pop removes entries. We record earlier pops below later pops in the heap. Thus the top of an arrow stack is read first, but the bottom of the recorded heap comes first.

There may be several exposed cycles to choose from. The point of Wilson's cycle-popping lemma is that these choices do not affect the answer. We use the following finite, deterministic form \cite{Wilson,PW,LyonsPeres}.

\begin{lemma}[Cycle popping]\label{lem:popping}
All complete popping sequences from a stack configuration $\sigma$ end at the same configuration $\sigma^\circ$ and record the same heap $H(\sigma)$. Moreover,
\[
\sigma\longmapsto(H(\sigma),\sigma^\circ)
\]
is a bijection from configurations of height $T$ to pairs $(H,\sigma^\circ)$ in which $\cont(H)\le T$ and $\sigma^\circ$ has height $T-\cont(H)$ and acyclic exposed function. It satisfies
\[
\wt(\sigma)=c^H\wt(\sigma^\circ).
\]
\end{lemma}

\begin{proof}
The key observation is that an exposed cycle persists until it is popped. Suppose $\delta$ is exposed initially, and consider any complete popping sequence. Some pop must remove an entry of $\delta$, or it would still be exposed at the end. At the first such pop, all its entries are still on top, so following the exposed arrows gives exactly $\delta$. Every earlier pop is disjoint from it. We may therefore move $\delta$ to the beginning without changing the heap or the final configuration. Comparing two complete sequences in this way and inducting on the number of entries proves order independence.

For the inverse, start with $(H,\sigma^\circ)$. At vertex $v$, put the successors of $v$ in the pieces of $H$ through $v$ above $\sigma^\circ_v$, in bottom-to-top heap order: the lowest piece supplies the top entry of the stack. Those pieces are totally ordered, so this is well defined. Any representative word of $H$ is now a legal popping sequence. It ends at $\sigma^\circ$, where no further pop is possible. The preceding paragraph shows that this is the inverse, and \eqref{eq:telescope} gives the weight identity.
\end{proof}

\begin{example}\label{ex:pop}
Take $\sigma_1=(2,*)$ and $\sigma_2=(1,3)$, top first, with all other stacks empty. The exposed arrows form $\alpha=(12)$. After it pops, the choices are $1\mapsto*$ and $2\mapsto3$, an acyclic function on $12$. Thus
\[
H(\sigma)=\alpha,\qquad
\wt(\sigma)=f_{12}\frac{A_1}{X_1}\frac{c_{23}X_3}{X_2}.
\]
Changing $\sigma_1$ to $(*,2)$ makes the initial exposed function acyclic. Stack order matters; the order of legal pops does not.
\end{example}

\subsection{Pop, peel, and repeat}
Once all exposed cycles have been popped, the top entries form an acyclic function, one of the terms in a principal minor. The entries below them need not be acyclic. We remove the top entries as one layer and repeat with the shorter stacks. We call this the \emph{pop--peel algorithm}.

Let $H_r$ be the heap popped in round $r$, and let $S_r$ be the support of the acyclic layer then removed. The sets $S_r$ decrease, so they are the layers of $U=\sum_r\one_{S_r}$. After peeling $S_{r-1}$, entries remain only at its vertices. This explains the one-step shift in the following definition.

\begin{definition}[Leveled heaps]\label{def:leveled}
A \emph{$U$-leveled heap} is a sequence $\bH=(H_1,H_2,\ldots)$ with
\[
H_r\in\Heap(U^{(r-1)})\qquad(r\ge1).
\]
Its weight and content are $c^{\bH}=\prod_rc^{H_r}$ and
$\cont(\bH)=\sum_r\cont(H_r)$.
\end{definition}

The heaps are eventually empty. In particular, $H_1$ is unrestricted, $H_2$ lives in $\supp U$, and $H_3$ lives in $U^{(2)}$. The level $r$ means that $r-1$ layers have already been peeled off before $H_r$ is popped; it is not the height of a piece within $H_r$. Notice also that $U$ records the layers \emph{left after popping}; it is not the initial stack height $T$.

\begin{theorem}[Leveled heap formula]\label{thm:A}
The pop--peel algorithm is a weight-preserving bijection
\[
\sigma\longleftrightarrow
\bigl(U,\bH,(\varphi_r)_{r\ge1}\bigr),
\]
where $\sigma$ has height $T$, $\bH$ is $U$-leveled with
$\cont(\bH)=T-U$, and $\varphi_r\in\Acyc(U^{(r)})$.
The weight on the right is $c^{\bH}\prod_r\wt(\varphi_r)$. Hence
\[
Y^T=\sum_U K_{T,U}B_U,
\qquad
K_{T,U}=\sum_{\substack{\bH\text{ $U$-leveled}\\\cont(\bH)=T-U}}c^{\bH}\in\NN[c].
\]
\end{theorem}

\begin{proof}
Let $T^{[r]}$ be the heights after round $r$, with $T^{[0]}=T$. Then
\[
T^{[r]}=T^{[r-1]}-\cont(H_r)-\one_{S_r},
\qquad
S_r=\supp\bigl(T^{[r-1]}-\cont(H_r)\bigr).
\]
Thus $S_{r+1}\subseteq\supp T^{[r]}\subseteq S_r$. The sets $S_r$ are the layers of $U$, and $T=U+\sum_r\cont(H_r)$. For $r\ge2$, the heap $H_r$ uses only vertices of $S_{r-1}$, so the output has the required form.

To reverse the algorithm, rebuild the stacks from the last round. Suppose the remaining stacks have support contained in $U^{(r)}$. Place $\varphi_r$ on top; the exposed function is now acyclic, on exactly $U^{(r)}$. Restore $H_r$ using Lemma~\ref{lem:popping}. The rebuilt stacks have support contained in $U^{(r-1)}$, so the construction continues. This reverses every round and recovers the initial height $T$.

The weights follow from the popping lemma. Summing over $\varphi_r$ gives $Y_{U^{(r)}}$ by Proposition~\ref{prop:forest}; their product is $B_U$.
\end{proof}

\begin{example}[A complete run]\label{ex:run}
Start with the following stacks, listed from top to bottom:
\[
\sigma_1=(*,2,*),\qquad
\sigma_2=(3,*,1,*),\qquad
\sigma_3=(2,4),\qquad\sigma_4=\emptyset.
\]
Their height is $T=(3,4,2,0)$. First the exposed arrows at $2$ and $3$ form $\beta$, which pops. The exposed function is then $\varphi_1=(1\mapsto*,2\mapsto*,3\mapsto4)$, acyclic on $123$ because the last arrow exits. Peeling this layer reveals $\alpha$. After it pops, the two roots on $12$ form $\varphi_2$. Thus
\[
(H_1,H_2)=(\beta,\alpha),\qquad
U=\one_{123}+\one_{12}=(2,2,1,0).
\]
The weight is $f_{23}f_{12}\wt(\varphi_1)\wt(\varphi_2)$. Figure~\ref{fig:run} also shows how to reverse the run.

\par\medskip
\begin{center}\begin{minipage}{\linewidth}\centering
\begin{tikzpicture}[baseline=(base)]
  \coordinate (base) at (0,0.84);
  \scell{1}{4}{black!24}{black!55}{$*$}
  \scell{1}{3}{calpha!22}{calpha!60!black}{$2$}
  \scell{1}{2}{black!8}{black!35}{$*$}
  \scell{2}{4}{cbeta!22}{cbeta!60!black}{$3$}
  \scell{2}{3}{black!24}{black!55}{$*$}
  \scell{2}{2}{calpha!22}{calpha!60!black}{$1$}
  \scell{2}{1}{black!8}{black!35}{$*$}
  \scell{3}{4}{cbeta!22}{cbeta!60!black}{$2$}
  \scell{3}{3}{black!24}{black!55}{$4$}
  \scols{4}{-0.2}
  \node[font=\tiny,text=black!55,anchor=east] at (-0.05,1.47) {top};
  \node[font=\scriptsize] at (1.0,-0.62) {stacks $\sigma$};
\end{tikzpicture}
\hspace{10pt}
\begin{tikzpicture}[baseline=(base)]
  \coordinate (base) at (0,0);
  \draw[-{Stealth[length=2mm]},black!60,line width=.6pt] (0,0.18) -- (1.7,0.18) node[midway,above,font=\scriptsize] {pop \& peel};
  \draw[{Stealth[length=2mm]}-,black!60,line width=.6pt] (0,-0.18) -- (1.7,-0.18) node[midway,below,font=\scriptsize] {restack};
\end{tikzpicture}
\hspace{10pt}
\begin{tikzpicture}[baseline=(base)]
  \coordinate (base) at (0,0.46);
  \fill[black!6] (0,0.46) rectangle (2.4,0.92);
  \heapgrid{4}{2}
  \heappiece{2}{3}{1}{cbeta}{\beta}
  \heappiece{1}{2}{2}{calpha}{\alpha}
  \node[font=\tiny,text=black!60,anchor=west] at (2.45,0.23) {$H_1$};
  \node[font=\tiny,text=black!60,anchor=west] at (2.45,0.69) {$H_2$};
  \node[font=\scriptsize] at (1.2,-0.62) {leveled heap};
\end{tikzpicture}
\hspace{14pt}
\begin{tikzpicture}[baseline=(base)]
  \coordinate (base) at (0,0.42);
  \scell{1}{1}{black!24}{black!55}{$*$}
  \scell{2}{1}{black!24}{black!55}{$*$}
  \scell{3}{1}{black!24}{black!55}{$4$}
  \scell{1}{2}{black!8}{black!35}{$*$}
  \scell{2}{2}{black!8}{black!35}{$*$}
  \node[font=\tiny,text=black!60,anchor=west] at (1.55,0.21) {$\varphi_1$};
  \node[font=\tiny,text=black!60,anchor=west] at (1.05,0.63) {$\varphi_2$};
  \scols{4}{-0.2}
  \node[font=\scriptsize] at (1.0,-0.62) {layers, $U=(2,2,1,0)$};
\end{tikzpicture}

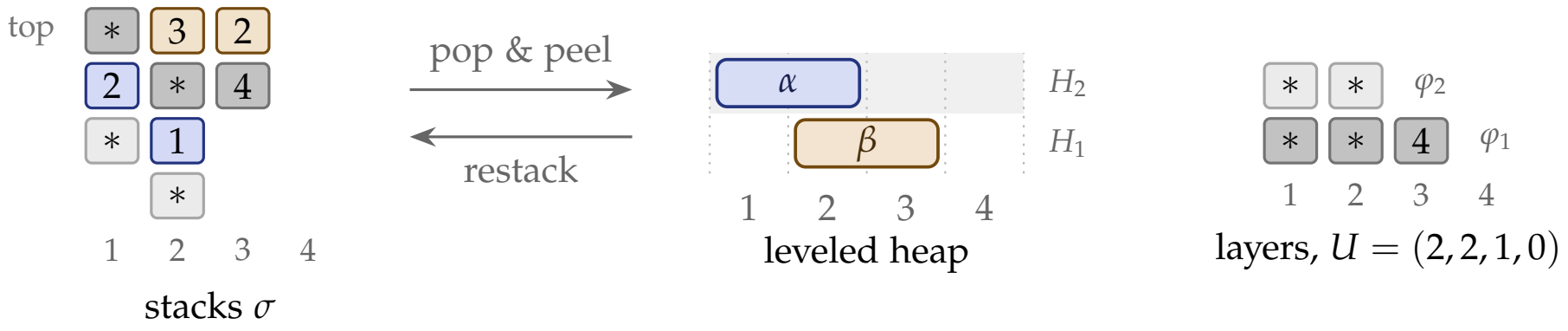
\captionof{figure}{The pop--peel bijection of Theorem~\ref{thm:A}. Left: stacks of height $(3,4,2,0)$, read top-down; a cell labeled $j$ is the arrow $v\to j$, and $*$ is a root. Cells are colored by their fate: $\beta$ pops first, then the dark layer $\varphi_1$ is peeled, then $\alpha$ pops, then the light layer $\varphi_2$ is peeled. Right: the output, read bottom-up: the $U$-leveled heap $(H_1,H_2)=(\beta,\alpha)$ and the layers in the diagram of $U$. Restacking---each level-$r$ heap placed above $\varphi_r$, and level $r$ above level $r+1$---recovers the stacks.}\label{fig:run}
\end{minipage}\end{center}
\medskip

\end{example}

\begin{example}[Where levels may occur]\label{ex:leveled}
For $U=(3,3,2,1)$, the nonempty layers are $1234$, $123$, and $12$. Thus the first two heaps may use all four vertices, the third may use only $123$, and the fourth may use only $12$. Figure~\ref{fig:leveled} shows
\[
(H_1,H_2,H_3,H_4)
=(\alpha\beta\gamma,\ \gamma\beta,\ \tau\alpha,\ \alpha^2).
\]
The two copies of $\alpha$ at level $4$ are allowed because that level uses $U^{(3)}=12$. The corresponding initial height is $T=U+\cont(\bH)$.

\par\medskip
\begin{center}\begin{minipage}{\linewidth}\centering
\begin{tikzpicture}[baseline=0pt]
  \ccellf{1}{1}{black!10}{black!40}{} \ccellf{2}{1}{black!10}{black!40}{} \ccellf{3}{1}{black!10}{black!40}{} \ccellf{4}{1}{black!10}{black!40}{}
  \ccellf{1}{2}{black!18}{black!45}{} \ccellf{2}{2}{black!18}{black!45}{} \ccellf{3}{2}{black!18}{black!45}{}
  \ccellf{1}{3}{black!27}{black!50}{} \ccellf{2}{3}{black!27}{black!50}{}
  \ccols{4}
  \node[font=\scriptsize] at (0.84,1.62) {$U=(3,3,2,1)$};
  \node[font=\tiny,anchor=west,text=black!70] at (1.76,0.21) {$U^{(1)}=\{1,2,3,4\}$};
  \node[font=\tiny,anchor=west,text=black!70] at (1.76,0.63) {$U^{(2)}=\{1,2,3\}$};
  \node[font=\tiny,anchor=west,text=black!70] at (1.76,1.05) {$U^{(3)}=\{1,2\}$};
\end{tikzpicture}\hspace{26pt}
\begin{tikzpicture}[baseline=0pt]
  \fill[black!6] (0,1.38) rectangle (2.4,2.30);
  \fill[black!6] (0,3.22) rectangle (2.4,4.14);
  \heapgrid{4}{9}
  \heappiece{1}{2}{1}{calpha}{\alpha}
  \heappiece{2}{3}{2}{cbeta}{\beta}
  \heappiece{3}{4}{3}{cgamma}{\gamma}
  \heappiece{3}{4}{4}{cgamma}{\gamma}
  \heappiece{2}{3}{5}{cbeta}{\beta}
  \heappiece{1}{3}{6}{ctau}{\tau}
  \heappiece{1}{2}{7}{calpha}{\alpha}
  \heappiece{1}{2}{8}{calpha}{\alpha}
  \heappiece{1}{2}{9}{calpha}{\alpha}
  \draw[decorate,decoration={brace,amplitude=3.5pt},black!55,line width=.5pt] (2.55,1.34) -- (2.55,0.04)
     node[midway,right=6pt,font=\scriptsize,text=black!75] {$H_1=\alpha\beta\gamma\in\Heap(U^{(0)})$};
  \draw[decorate,decoration={brace,amplitude=3.5pt},black!55,line width=.5pt] (2.55,2.26) -- (2.55,1.42)
     node[midway,right=6pt,font=\scriptsize,text=black!75] {$H_2=\gamma\beta\in\Heap(U^{(1)})$};
  \draw[decorate,decoration={brace,amplitude=3.5pt},black!55,line width=.5pt] (2.55,3.18) -- (2.55,2.34)
     node[midway,right=6pt,font=\scriptsize,text=black!75] {$H_3=\tau\alpha\in\Heap(U^{(2)})$};
  \draw[decorate,decoration={brace,amplitude=3.5pt},black!55,line width=.5pt] (2.55,4.10) -- (2.55,3.26)
     node[midway,right=6pt,font=\scriptsize,text=black!75] {$H_4=\alpha^2\in\Heap(U^{(3)})$};
\end{tikzpicture}

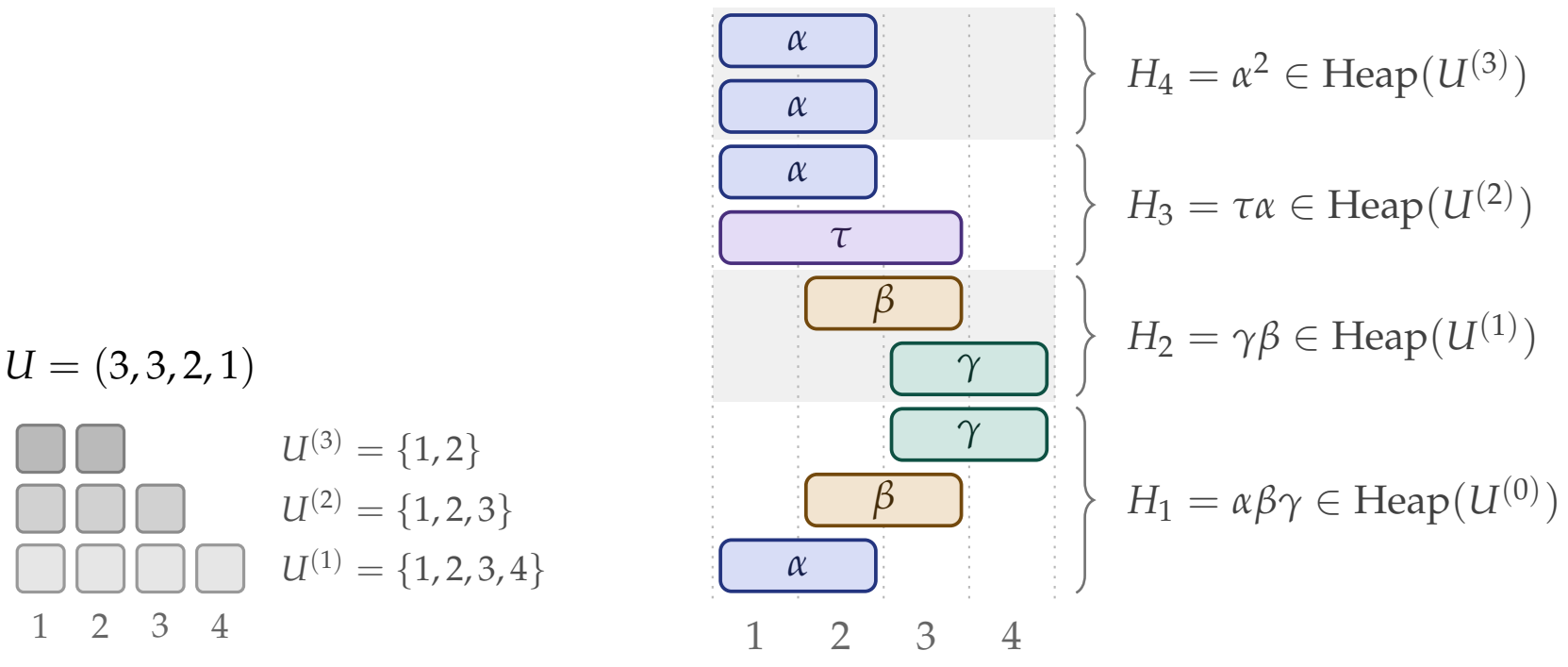
\captionof{figure}{Left: the diagram of $U=(3,3,2,1)$, whose rows are its layers. Right: a $U$-leveled heap drawn as one stack of bars, with alternating bands marking the levels.}\label{fig:leveled}
\end{minipage}\end{center}
\medskip

\end{example}

\begin{example}[A coefficient of eight]\label{ex:eight}
Let us count all the ways to obtain $U=(2,2,1,0)$ from $T=(3,4,2,0)$. We have $T-U=\cont(\alpha)+\cont(\beta)$. The cycles are fixed, but their levels are not: $\alpha$ may occur at levels $1,2,3$, while $\beta$ may occur only at levels $1,2$.

There are four same-level possibilities: both pieces lie at level $1$ or $2$, in either heap order. There are four different-level possibilities, with
\[
(\text{level of }\alpha,\text{level of }\beta)
=(1,2),(2,1),(3,1),(3,2).
\]
These are all the possibilities, since no other collection of cycle supports has content $(1,2,1,0)$. Consequently,
\[
K_{(3,4,2,0),(2,2,1,0)}=8f_{12}f_{23}.
\]
The run in Example~\ref{ex:run} is the sixth configuration in Figure~\ref{fig:eight}.

\par\medskip
\begin{center}
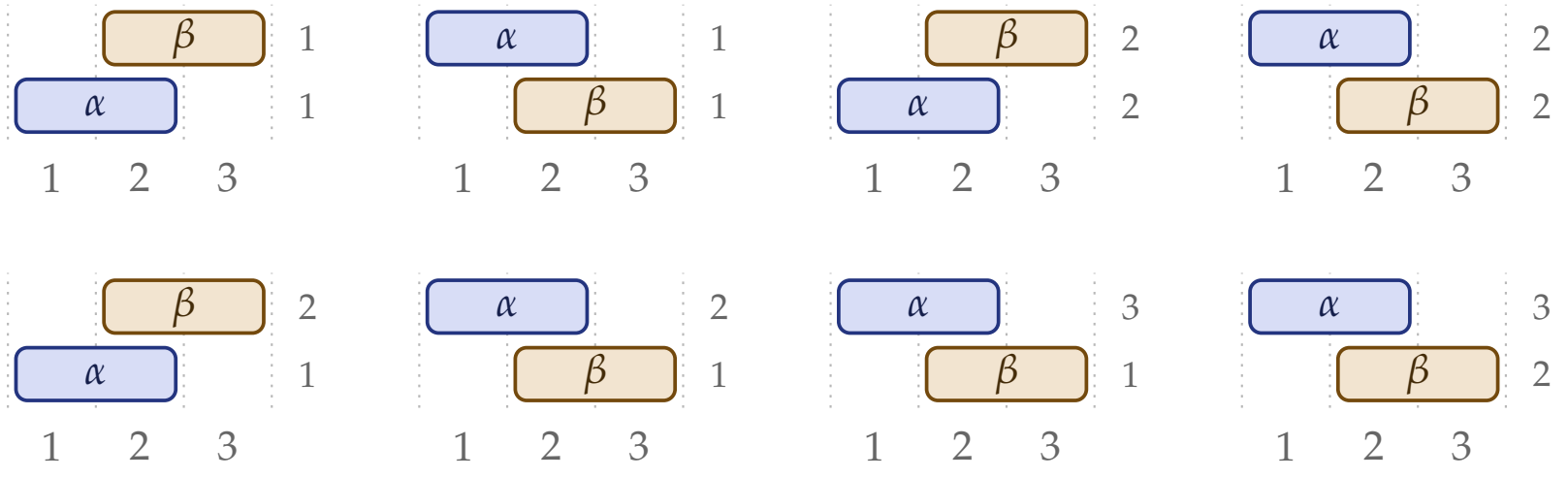
\begin{minipage}{\linewidth}\centering
\newcommand\eightAB[2]{\begin{tikzpicture}[baseline=0pt]
  \heapgrid{3}{2}\heappiece{1}{2}{1}{calpha}{\alpha}\heappiece{2}{3}{2}{cbeta}{\beta}\lvl{3}{1}{#1}\lvl{3}{2}{#2}
\end{tikzpicture}}%
\newcommand\eightBA[2]{\begin{tikzpicture}[baseline=0pt]
  \heapgrid{3}{2}\heappiece{2}{3}{1}{cbeta}{\beta}\heappiece{1}{2}{2}{calpha}{\alpha}\lvl{3}{1}{#1}\lvl{3}{2}{#2}
\end{tikzpicture}}%
\eightAB{1}{1}\hspace{16pt}\eightBA{1}{1}\hspace{16pt}\eightAB{2}{2}\hspace{16pt}\eightBA{2}{2}\\[10pt]
\eightAB{1}{2}\hspace{16pt}\eightBA{1}{2}\hspace{16pt}\eightBA{1}{3}\hspace{16pt}\eightBA{2}{3}
\captionof{figure}{The eight $(2,2,1,0)$-leveled heaps of content $\cont(\alpha)+\cont(\beta)$; the number beside each piece is its level. Pieces at the same level form a heap $\alpha\beta$ or $\beta\alpha$; pieces at different levels are ordered by level.}\label{fig:eight}
\end{minipage}\end{center}
\medskip

\end{example}

\section{The pair rule}\label{sec:pair}

The leveled-heap formula starts with arbitrary stacks. In a product $Y_IY_J$, each factor already forbids cycles inside its own set. We now derive a rule that remembers these restrictions: one becomes a condition at the bottom of a heap, and the other a condition at the top.

\Needspace{7\baselineskip}
\subsection{Boundary heaps}
\begin{definition}\label{def:boundary}
Let $\kappa=\one_I+\one_J$. An \emph{$(I,J)$-boundary heap} is a heap $H$ with $\cont(H)\le\kappa$ such that
\[
\max(H)\cap\cC(I)=\emptyset,
\qquad
\min(H)\cap\cC(J)=\emptyset.
\]
Here intersections mean that we inspect the cycle labeling each extremal occurrence.
\end{definition}

Every piece of $H$ is supported in $I\cup J$. The residual vector
$V=\kappa-\cont(H)$ has entries at most two, so
\[
B_V=Y_{V^{(1)}}Y_{V^{(2)}}.
\]
If $I\subseteq J$ or $J\subseteq I$, only the empty heap is possible, as expected for a product that is already nested. Reversing a heap exchanges the boundary conditions, so the rule is symmetric in $I,J$ despite the asymmetric indexing.

\begin{theorem}[Pair rule]\label{thm:B}
For all $I,J\subseteq[n]$,
\[
\boxed{
Y_IY_J=\sum_{H\text{ an $(I,J)$-boundary heap}}
 c^H B_{\one_I+\one_J-\cont(H)}.
}
\]
\end{theorem}

\begin{example}\label{ex:pair}
For $I=12$ and $J=23$, the available multiplicities are $(1,2,1,0)$. The one-piece heaps $\alpha$ and $\beta$ fail the top and bottom conditions, respectively. The heap $\beta\alpha$ fails both. The permitted heaps are
\[
\emptyset,\qquad\tau,\qquad\alpha\beta,
\]
and they give
\[
Y_{12}Y_{23}=Y_{123}Y_2+tY_2+f_{12}f_{23}.
\]
For the reversed pair, the two-piece heap is $\beta\alpha$ instead. Figure~\ref{fig:pair} shows the residual layers.

\par\medskip
\begin{center}
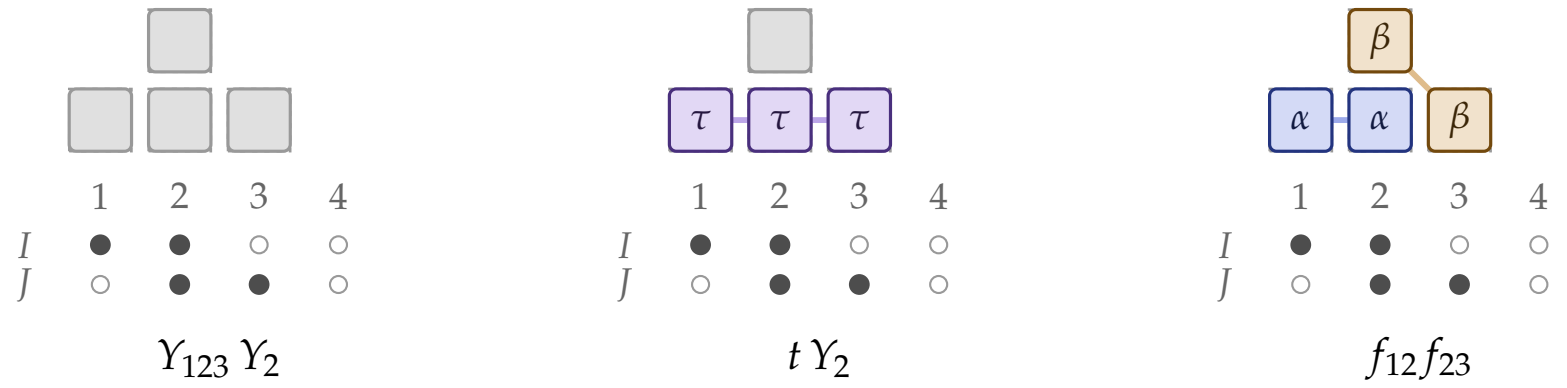
\begin{minipage}{\linewidth}\centering
\begin{tikzpicture}[baseline=0pt,scale=1.25]
  \cslot{1}{1}\cslot{2}{1}\cslot{2}{2}\cslot{3}{1}
  \ccellf{1}{1}{black!12}{black!40}{}\ccellf{2}{1}{black!12}{black!40}{}\ccellf{2}{2}{black!12}{black!40}{}\ccellf{3}{1}{black!12}{black!40}{}
  \ccols{4}
  \membrow{0.42}{4}{-0.45}{$I$}{1,2}\membrow{0.42}{4}{-0.66}{$J$}{2,3}
  \node[font=\scriptsize] at (0.84,-1.05) {$Y_{123}\,Y_2$};
\end{tikzpicture}\hspace{40pt}
\begin{tikzpicture}[baseline=0pt,scale=1.25]
  \cslot{1}{1}\cslot{2}{1}\cslot{2}{2}\cslot{3}{1}
  \cjoin{1}{1}{3}{1}{ctau}
  \ccell{1}{1}{ctau}{\tau}\ccell{2}{1}{ctau}{\tau}\ccell{3}{1}{ctau}{\tau}
  \ccellf{2}{2}{black!12}{black!40}{}
  \ccols{4}
  \membrow{0.42}{4}{-0.45}{$I$}{1,2}\membrow{0.42}{4}{-0.66}{$J$}{2,3}
  \node[font=\scriptsize] at (0.84,-1.05) {$t\,Y_2$};
\end{tikzpicture}\hspace{40pt}
\begin{tikzpicture}[baseline=0pt,scale=1.25]
  \cslot{1}{1}\cslot{2}{1}\cslot{2}{2}\cslot{3}{1}
  \cjoin{1}{1}{2}{1}{calpha}\cjoin{2}{2}{3}{1}{cbeta}
  \ccell{1}{1}{calpha}{\alpha}\ccell{2}{1}{calpha}{\alpha}
  \ccell{2}{2}{cbeta}{\beta}\ccell{3}{1}{cbeta}{\beta}
  \ccols{4}
  \membrow{0.42}{4}{-0.45}{$I$}{1,2}\membrow{0.42}{4}{-0.66}{$J$}{2,3}
  \node[font=\scriptsize] at (0.84,-1.05) {$f_{12}f_{23}$};
\end{tikzpicture}
\captionof{figure}{The $(\{1,2\},\{2,3\})$-boundary heaps. Dashed boxes are the capacity $\kappa=\one_I+\one_J$; each piece fills one box in each of its columns, in heap order (cells of one piece are joined), and the dots below record membership in $I$ and $J$. The unused boxes (gray) form $V=\kappa-\cont(H)$, read as $Y_{V^{(1)}}Y_{V^{(2)}}$.}\label{fig:pair}
\end{minipage}\end{center}
\medskip

\end{example}

\subsection{Two levels are enough}
We will expand $Y^{\one_I+\one_J}$ in two ways. Only the capacity-two case of Theorem~\ref{thm:A} is needed: every stack has at most two entries, so the algorithm has at most two popping rounds.

\begin{lemma}\label{lem:twolevels}
Let $\bH$ be $U$-leveled with $\cont(\bH)=\kappa-U$, where $\kappa=\one_I+\one_J$. Then $H_r=\emptyset$ for $r\ge3$. The second level $H_2$ is a clique supported in
\[
(I\cap J)\setminus\supp\cont(H_1).
\]
Conversely, any heap $H_1$ with $\cont(H_1)\le\kappa$, together with such a clique $H_2$, is a leveled heap for
$U=\kappa-\cont(H_1)-\cont(H_2)$.
\end{lemma}

\begin{proof}
A piece at level $r$ through $v$ requires $U_v\ge r-1$ and consumes one additional unit of $\kappa_v\le2$. Thus $r\le2$. At level $2$, this forces $U_v=1$ and $\kappa_v=2$, leaving no room for any other piece through $v$. The converse follows from the same count.
\end{proof}

Thus $H_2$ uses only shared vertices untouched by $H_1$. We obtain
\begin{equation}\label{eq:Ykappa}
Y^\kappa=\sum_{(H_1,H_2)}c^{H_1}c^{H_2}
 B_{\kappa-\cont(H_1)-\cont(H_2)},
\end{equation}
with the sum over the pairs in Lemma~\ref{lem:twolevels}.

\subsection{Stripping the boundaries}
To obtain a boundary heap, we would like to remove the minimal pieces lying in $J$ and the maximal pieces lying in $I$. Removing them may expose further pieces of these kinds. The useful observation is that we do not have to keep removing pieces: we shrink $I$ and $J$ by the vertices already used. A newly exposed piece then meets one of the deleted vertices, so it is no longer forbidden. The following lemma makes this precise.

\begin{lemma}[Boundary stripping]\label{lem:strip}
The pairs $(H_1,H_2)$ in \eqref{eq:Ykappa} are in weight- and
content-preserving bijection with triples $(K_I,K_J,H)$, where
$K_I\subseteq\cC(I)$ and $K_J\subseteq\cC(J)$ are cliques and $H$
is an $(I',J')$-boundary heap for
$I'=I\setminus V(K_I)$ and $J'=J\setminus V(K_J)$.
In particular,
\begin{equation}\label{eq:strip-content}
\kappa-\cont(H_1)-\cont(H_2)
=\one_{I'}+\one_{J'}-\cont(H).
\end{equation}
\end{lemma}

\begin{proof}
The forward map removes pieces from the two ends of $H_1$:
\begin{enumerate}[leftmargin=2.2em,label=\textup{(\roman*)}]
\item Let $K_J$ be the original minimal occurrences of $H_1$ supported in $J$.
\item Let $E_I$ be the original maximal occurrences of $H_1$ supported in $I$ that were not assigned to $K_J$. Set $K_I=E_I\sqcup H_2$.
\item Delete $K_J$ and $E_I$ from $H_1$, leaving $H$.
\end{enumerate}
To reverse the map, put $K_J$ below $H$ and $K_I$ above it, except for the pieces that belong at level $2$. These are the shared pieces of $K_I$ disjoint from everything else. Put
\[
Q=\{\gamma\in K_I:V(\gamma)\subseteq J,\quad
V(\gamma)\cap(V(K_J)\cup\supp\cont(H))=\emptyset\}.
\]
Then
\begin{equation}\label{eq:strip-inverse}
H_2=Q,\qquad H_1=K_J\cdot H\cdot(K_I\setminus Q).
\end{equation}

Let us check the construction. The minimal pieces form a clique, as do the maximal pieces, and $H_2$ is disjoint from $H_1$ by Lemma~\ref{lem:twolevels}. Hence $K_I$ and $K_J$ are cliques. The removed occurrences can be placed at the two ends of a representative word, so the remainder is again a heap. Counting the removed occurrences gives weight preservation and \eqref{eq:strip-content}; in particular, $\cont(H)\le\one_{I'}+\one_{J'}$.

Consider a minimal piece $p$ of $H$. If it was already minimal in $H_1$, it was not removed, so its support is not contained in $J$. Otherwise an immediate predecessor $q$ in $H_1$ was removed. Such a predecessor meets $p$ and cannot have been maximal, so $q\in K_J$. Thus $p$ uses a vertex deleted from $J$, and is not supported in $J'$. This proves the bottom condition. The top condition is identical: a newly exposed maximal piece meets a removed piece of $E_I$, so it is not supported in $I'$.

Now start with a triple $(K_I,K_J,H)$ and use \eqref{eq:strip-inverse}. The pieces of $Q$ lie in $I\cap J$ and are disjoint from the new $H_1$, so the two-level conditions hold. We claim that stripping again recovers the original cliques. Certainly every piece of $K_J$ is minimal and lies in $J$. Any other minimal piece avoids $V(K_J)$. If it comes from $H$, its support cannot lie in $J$, since it would then lie in $J'$; if it comes from $K_I$, being minimal and supported in $J$ would put it in $Q$. Thus we recover exactly $K_J$. The same argument at the top recovers $K_I\setminus Q$, and adding $H_2=Q$ recovers $K_I$.

Conversely, for a pair $(H_1,H_2)$, every piece of $H_2$ belongs to $Q$. A piece of $E_I$ cannot belong to $Q$: it would be supported in $J$ and disjoint from the rest of $H_1$, making it a minimal $J$-piece already assigned to $K_J$. Thus $Q=H_2$, and reattaching the removed pieces recovers $H_1$.
\end{proof}

An isolated shared occurrence in level $1$ goes to $K_J$; an occurrence in level $2$ goes to $K_I$. This convention is the only asymmetry. If the same cycle occurs twice, the occurrences must be treated separately.

\begin{example}\label{ex:strip}
Take $I=1234$, $J=123$, and $(H_1,H_2)=(\alpha\beta\gamma,\emptyset)$. Stripping gives
\[
K_J=\{\alpha\},\qquad H=\beta,\qquad K_I=\{\gamma\},
\]
with reduced sets $I'=12$, $J'=3$. The piece $\beta$ is now both minimal and maximal, but it lies in neither reduced set. This is where shrinking the sets matters. The residual vector is $(1,0,0,0)$ before and after stripping.

For the shared-cycle convention, take $I=J=12$. There are four pairs:
\[
(\emptyset,\emptyset),\quad(\alpha,\emptyset),\quad
(\emptyset,\alpha),\quad(\alpha^2,\emptyset).
\]
They assign no cycle, a cycle to $K_J$, a cycle to $K_I$, or one occurrence to each clique, respectively. In all four cases $H=\emptyset$. Figures~\ref{fig:strip} and~\ref{fig:shared} show the two cases.

\par\medskip
\begin{center}\begin{minipage}{\linewidth}\centering
\begin{tikzpicture}[baseline=(base)]
  \coordinate (base) at (0,0.5);
  \heapgrid{4}{3}
  \heappieceK{1}{2}{1}{calpha}{\alpha}
  \heappiece{2}{3}{2}{cbeta}{\beta}
  \heappieceKI{3}{4}{3}{cgamma}{\gamma}
  \node[font=\tiny,text=black!65,anchor=east] at (-0.05,0.23) {$K_J$};
  \node[font=\tiny,text=black!65,anchor=west] at (2.45,1.15) {$K_I$};
  \membrow{0.6}{4}{-0.52}{$I$}{1,2,3,4}\membrow{0.6}{4}{-0.74}{$J$}{1,2,3}
  \node[font=\scriptsize] at (1.2,-1.12) {$(H_1,H_2)=(\alpha\beta\gamma,\emptyset)$};
\end{tikzpicture}
\hspace{8pt}
\begin{tikzpicture}[baseline=(base)]
  \coordinate (base) at (0,0);
  \draw[-{Stealth[length=2mm]},black!60,line width=.6pt] (0,0.18) -- (1.6,0.18) node[midway,above,font=\scriptsize] {strip};
  \draw[{Stealth[length=2mm]}-,black!60,line width=.6pt] (0,-0.18) -- (1.6,-0.18) node[midway,below,font=\scriptsize,align=center] {reattach};
\end{tikzpicture}
\hspace{8pt}
\begin{tikzpicture}[baseline=(base)]
  \coordinate (base) at (0,0.5);
  \begin{scope}[yshift=1.05cm]
    \heappieceKI{3}{4}{1}{cgamma}{\gamma}
    \node[font=\tiny,text=black!65,anchor=west] at (2.45,0.23) {$K_I=\{\gamma\}$};
  \end{scope}
  \heapgrid{4}{1}
  \heappiece{2}{3}{1}{cbeta}{\beta}
  \node[font=\tiny,text=black!65,anchor=west] at (2.45,0.23) {$H=\beta$};
  \membrow{0.6}{4}{-0.52}{$I'$}{1,2}\membrow{0.6}{4}{-0.74}{$J'$}{3}
  \begin{scope}[yshift=-1.42cm]
    \heappieceK{1}{2}{1}{calpha}{\alpha}
    \node[font=\tiny,text=black!65,anchor=west] at (2.45,0.23) {$K_J=\{\alpha\}$};
  \end{scope}
\end{tikzpicture}


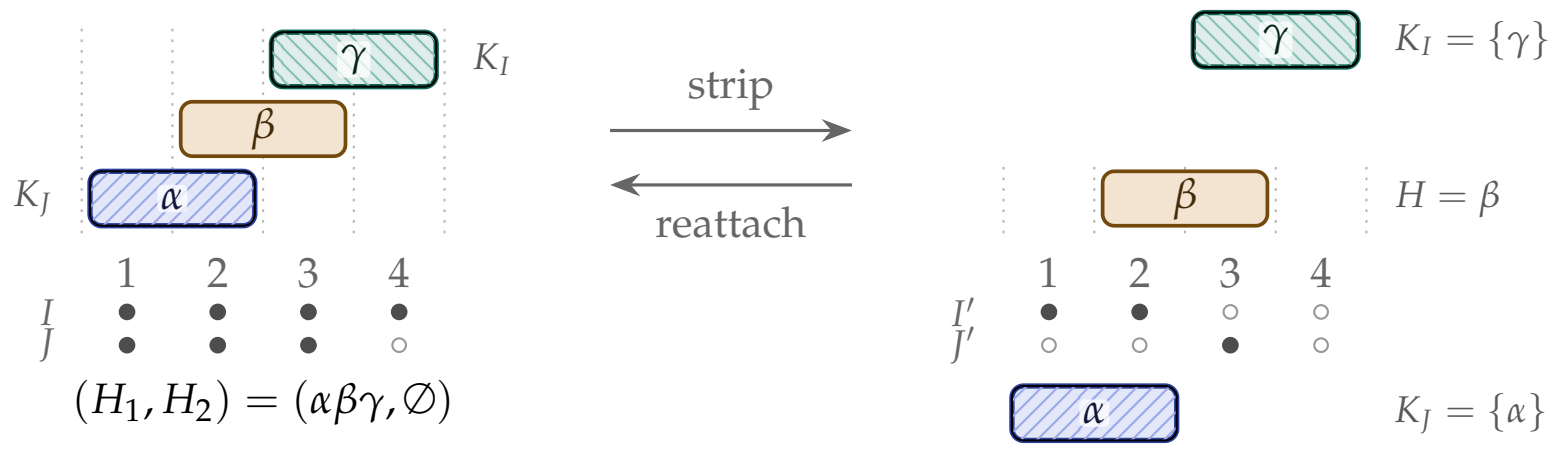
\captionof{figure}{Boundary stripping in the example $I=1234$, $J=123$. The original minimal $J$-piece $\alpha$ goes to $K_J$ and the maximal $I$-piece $\gamma$ goes to $K_I$; the remaining $\beta$ is a boundary heap for the reduced sets. The inverse reattaches the two cliques.}\label{fig:strip}
\end{minipage}\end{center}
\medskip

\begin{center}
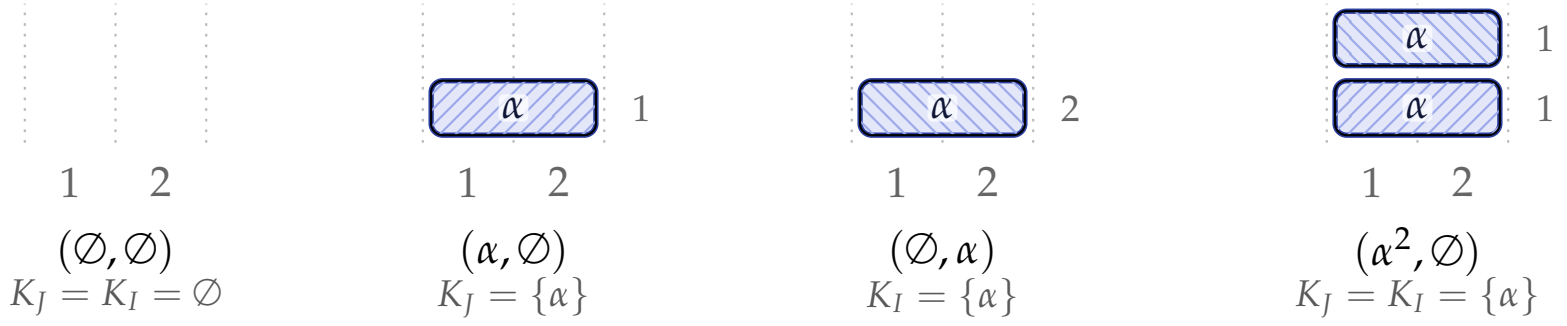
\begin{minipage}{\linewidth}\centering

\newcommand\sharedpic[3]{\begin{tikzpicture}[baseline=0pt]
  \heapgrid{2}{2}#1
  \node[font=\scriptsize] at (0.6,-0.72) {#2};
  \node[font=\tiny,text=black!65] at (0.6,-1.02) {#3};
\end{tikzpicture}}%
\sharedpic{}{$(\emptyset,\emptyset)$}{$K_J=K_I=\emptyset$}\hspace{30pt}
\sharedpic{\heappieceK{1}{2}{1}{calpha}{\alpha}\lvl{2}{1}{1}}{$(\alpha,\emptyset)$}{$K_J=\{\alpha\}$}\hspace{30pt}
\sharedpic{\heappieceKI{1}{2}{1}{calpha}{\alpha}\lvl{2}{1}{2}}{$(\emptyset,\alpha)$}{$K_I=\{\alpha\}$}\hspace{30pt}
\sharedpic{\heappieceK{1}{2}{1}{calpha}{\alpha}\heappieceKI{1}{2}{2}{calpha}{\alpha}\lvl{2}{1}{1}\lvl{2}{2}{1}}{$(\alpha^2,\emptyset)$}{$K_J=K_I=\{\alpha\}$}

\captionof{figure}{The shared-cycle convention for $I=J=12$. A level-1 isolated piece enters $K_J$, and a level-2 piece enters $K_I$; when two occurrences are present, the lower goes to $K_J$ and the upper to $K_I$.}\label{fig:shared}
\end{minipage}\end{center}
\medskip

\end{example}

\subsection{Two expansions of the same monomial}
We now compare removing cycles separately from the two functions with popping all the stacks together.

\begin{proof}[Proof of Theorem~\ref{thm:B}]
Let $R(I,J)$ denote the proposed boundary-heap sum. We induct on $|I|+|J|$, starting with the empty pair. The cycle decomposition \eqref{eq:cycdec}, applied to both factors, gives
\[
Y^\kappa=\sum_{K_I,K_J}c^{K_I}c^{K_J}Y_{I'}Y_{J'}.
\]
On the other hand, Theorem~\ref{thm:A} and boundary stripping give
\[
Y^\kappa=\sum_{K_I,K_J}c^{K_I}c^{K_J}R(I',J').
\]
The sums use the same cliques and reduced sets. If either clique is nonempty, the corresponding summands agree by induction, since $|I'|+|J'|<|I|+|J|$. Canceling them leaves the terms with both cliques empty, whose equality is $Y_IY_J=R(I,J)$.
\end{proof}

For $I=J=12$, the four terms are simply
\[
Y_1^2Y_2^2=(Y_{12}+f_{12})^2
=Y_{12}^2+2f_{12}Y_{12}+f_{12}^2.
\]
The two placements of one shared cycle explain the coefficient $2$.

\section{Positive straightening}\label{sec:straightening}

\subsection{Products of principal minors}
We now apply the pair rule repeatedly. A nonempty heap lowers total vertex multiplicity. The empty heap does not, but it makes the sets more nested by replacing them with their union and intersection. The proof measures these two kinds of progress to show that straightening terminates.

\begin{theorem}[Positive straightening]\label{thm:positivity}
Every product of weighted principal minors is a linear combination of the $B_V$ with coefficients in $\NN[c]$. In particular, all multiplication coefficients in the basis $\{B_U\}$ are nonnegative.
\end{theorem}

\begin{proof}
Consider $\prod_aY_{I_a}$. If its index sets form a chain, it is already $B_U$ for $U=\sum_a\one_{I_a}$. Otherwise apply Theorem~\ref{thm:B} to two incomparable sets $I,J$.

Track
\[
N=\sum_a|I_a|,\qquad Q=\sum_a|I_a|^2.
\]
A nonempty heap strictly decreases $N$. The empty heap replaces $I,J$ by $I\cup J,I\cap J$; it preserves $N$ and increases $Q$ by
\[
2|I\setminus J|\,|J\setminus I|>0.
\]
Since $0\le Q\le N^2$, there are only finitely many steps on each branch, with a uniform bound depending on the initial $N$. Every branch ends at a chain product, and every coefficient produced along it is nonnegative.
\end{proof}

\begin{example}\label{ex:straightening}
In the running graph,
\[
\begin{aligned}
Y_{12}Y_{23}Y_3
&=Y_{123}Y_2Y_3+tY_2Y_3+f_{12}f_{23}Y_3\\
&=Y_{123}Y_{23}+f_{23}Y_{123}
+tY_{23}+tf_{23}+f_{12}f_{23}Y_3.
\end{aligned}
\]
The first line uses \eqref{eq:intro}; the second uses $Y_2Y_3=Y_{23}+f_{23}$. All five final products are nested. In the empty-heap branch of the first step, the total size stays $5$, while $Q$ changes from $2^2+2^2+1^2=9$ to $3^2+1^2+1^2=11$.
\end{example}

Different choices of pairs can give different heap descriptions of the same coefficient. The basis makes the resulting coefficient unique, but we do not know a one-step heap rule for a general product $B_UB_W$.

\subsection{Exit paths}
So far we have handled products of $Y$-variables. To include the initial variables $X_i$, we need one more operation: following a path rather than popping a cycle.

A monomial $X^aB_U$ is \emph{compatible} if
\[
a_i>0\quad\Longrightarrow\quad U_i=0.
\]
At $c_{ij}=1$, these are exactly the cluster monomials of $\cA_\Gamma$ \cite{LP2,DTW}. To straighten an incompatible product, we follow an arrow path out of its outermost minor.

\begin{lemma}[Exit-path identity {\cite[Lemma~4.7]{LP2}}]\label{lem:exitpath}
For $i\in K\subseteq[n]$,
\begin{equation}\label{eq:exitpath}
\begin{aligned}
X_iY_K={}&\sum_{p:i\to k\text{ in }\Gamma[K]}
 c^pA_kY_{K\setminus V(p)}\\
&+\sum_{p:i\to k\text{ in }\Gamma[K]}
 \sum_{\substack{k\to j\\j\notin K}}
 c^pc_{kj}X_jY_{K\setminus V(p)}.
\end{aligned}
\end{equation}
Here $p$ is simple, the one-vertex path is allowed, and $c^p$ is the product of its edge weights.
\end{lemma}

\begin{proof}
Take an acyclic function on $K$ and follow its arrows from $i$. The resulting path is simple and ends either at a root $k$ or at an exit $k\to j$ with $j\notin K$. After multiplying by $X_i$, its weight telescopes to $c^pA_k$ or $c^pc_{kj}X_j$.

The remaining choices are an acyclic function on $K\setminus V(p)$. Conversely, any such function can be joined to the path and its terminal choice. No new cycle is created: a path entering $p$ must follow it to its root or exit. Sum over the remaining function using Proposition~\ref{prop:forest}.
\end{proof}

\begin{example}\label{ex:exitpath}
For $K=123$ and $i=1$, the possible paths are $(1)$, $1\to2$, and $1\to2\to3$. Only the last has an exit, along $3\to4$. Thus
\[
\boxed{
X_1Y_{123}=A_1Y_{23}+c_{12}A_2Y_3
+c_{12}c_{23}A_3+c_{12}c_{23}c_{34}X_4.
}
\]
Figure~\ref{fig:exitpath} shows the four choices.

\par\medskip
\begin{center}\begin{minipage}{\linewidth}\centering
\begin{tabular}{@{}c@{\hspace{28pt}}c@{}}
\begin{tikzpicture}[scale=0.95]
  \fill[inkblue!6,rounded corners=10pt] (-0.5,-1.0) rectangle (3.15,0.55);
  \draw[inkblue!35,rounded corners=10pt,line width=.5pt,dash pattern=on 2pt off 1.6pt] (-0.5,-1.0) rectangle (3.15,0.55);
  \node[font=\scriptsize,text=inkblue!75] at (-0.24,0.36) {$K$};
  \node[gvertsm] (u1) at (0.00,0) {$1$};
  \node[gvertsm] (u2) at (1.32,0) {$2$};
  \node[gvertsm] (u3) at (2.64,0) {$3$};
  \node[gvertsm] (u4) at (3.96,0) {$4$};
  \draw[gdim] (u1) to[bend left=20] (u2);
  \draw[gdim] (u2) to[bend left=20] (u1);
  \draw[gdim] (u2) to[bend left=20] (u3);
  \draw[gdim] (u3) to[bend left=20] (u2);
  \draw[gdim] (u3) to[bend left=20] (u4);
  \draw[gdim] (u4) to[bend left=20] (u3);
  \draw[gdim] (u3) to[bend left=40] (u1);
  \draw[inkblue!80,line width=.9pt] (u1.center) circle (0.235);
  \node[font=\small] at (1.98,-1.44) {$p=(1)\ \mapsto\ A_1\,Y_{23}$};
\end{tikzpicture} &
\begin{tikzpicture}[scale=0.95]
  \fill[inkblue!6,rounded corners=10pt] (-0.5,-1.0) rectangle (3.15,0.55);
  \draw[inkblue!35,rounded corners=10pt,line width=.5pt,dash pattern=on 2pt off 1.6pt] (-0.5,-1.0) rectangle (3.15,0.55);
  \node[font=\scriptsize,text=inkblue!75] at (-0.24,0.36) {$K$};
  \node[gvertsm] (u1) at (0.00,0) {$1$};
  \node[gvertsm] (u2) at (1.32,0) {$2$};
  \node[gvertsm] (u3) at (2.64,0) {$3$};
  \node[gvertsm] (u4) at (3.96,0) {$4$};
  \draw[gdim] (u1) to[bend left=20] (u2);
  \draw[gdim] (u2) to[bend left=20] (u1);
  \draw[gdim] (u2) to[bend left=20] (u3);
  \draw[gdim] (u3) to[bend left=20] (u2);
  \draw[gdim] (u3) to[bend left=20] (u4);
  \draw[gdim] (u4) to[bend left=20] (u3);
  \draw[gdim] (u3) to[bend left=40] (u1);
  \draw[ghl,draw=inkblue!80] (u1) to[bend left=20] (u2);
  \draw[inkblue!80,line width=.9pt] (u2.center) circle (0.235);
  \node[font=\small] at (1.98,-1.44) {$p=1\to2\ \mapsto\ c_{12}A_2\,Y_3$};
\end{tikzpicture} \\[12pt]
\begin{tikzpicture}[scale=0.95]
  \fill[inkblue!6,rounded corners=10pt] (-0.5,-1.0) rectangle (3.15,0.55);
  \draw[inkblue!35,rounded corners=10pt,line width=.5pt,dash pattern=on 2pt off 1.6pt] (-0.5,-1.0) rectangle (3.15,0.55);
  \node[font=\scriptsize,text=inkblue!75] at (-0.24,0.36) {$K$};
  \node[gvertsm] (u1) at (0.00,0) {$1$};
  \node[gvertsm] (u2) at (1.32,0) {$2$};
  \node[gvertsm] (u3) at (2.64,0) {$3$};
  \node[gvertsm] (u4) at (3.96,0) {$4$};
  \draw[gdim] (u1) to[bend left=20] (u2);
  \draw[gdim] (u2) to[bend left=20] (u1);
  \draw[gdim] (u2) to[bend left=20] (u3);
  \draw[gdim] (u3) to[bend left=20] (u2);
  \draw[gdim] (u3) to[bend left=20] (u4);
  \draw[gdim] (u4) to[bend left=20] (u3);
  \draw[gdim] (u3) to[bend left=40] (u1);
  \draw[ghl,draw=inkblue!80] (u1) to[bend left=20] (u2);
  \draw[ghl,draw=inkblue!80] (u2) to[bend left=20] (u3);
  \draw[inkblue!80,line width=.9pt] (u3.center) circle (0.235);
  \node[font=\small] at (1.98,-1.44) {$p=1\to2\to3\ \mapsto\ c_{12}c_{23}A_3$};
\end{tikzpicture} &
\begin{tikzpicture}[scale=0.95]
  \fill[inkblue!6,rounded corners=10pt] (-0.5,-1.0) rectangle (3.15,0.55);
  \draw[inkblue!35,rounded corners=10pt,line width=.5pt,dash pattern=on 2pt off 1.6pt] (-0.5,-1.0) rectangle (3.15,0.55);
  \node[font=\scriptsize,text=inkblue!75] at (-0.24,0.36) {$K$};
  \node[gvertsm] (u1) at (0.00,0) {$1$};
  \node[gvertsm] (u2) at (1.32,0) {$2$};
  \node[gvertsm] (u3) at (2.64,0) {$3$};
  \node[gvertsm] (u4) at (3.96,0) {$4$};
  \draw[gdim] (u1) to[bend left=20] (u2);
  \draw[gdim] (u2) to[bend left=20] (u1);
  \draw[gdim] (u2) to[bend left=20] (u3);
  \draw[gdim] (u3) to[bend left=20] (u2);
  \draw[gdim] (u3) to[bend left=20] (u4);
  \draw[gdim] (u4) to[bend left=20] (u3);
  \draw[gdim] (u3) to[bend left=40] (u1);
  \draw[ghl,draw=inkblue!80] (u1) to[bend left=20] (u2);
  \draw[ghl,draw=inkblue!80] (u2) to[bend left=20] (u3);
  \draw[ghl,draw=cbeta] (u3) to[bend left=20] (u4);
  \draw[inkblue!80,line width=.9pt] (u3.center) circle (0.235);
  \node[font=\small] at (1.98,-1.44) {exit $3\to4\ \mapsto\ c_{12}c_{23}c_{34}\,X_4$};
\end{tikzpicture}
\end{tabular}

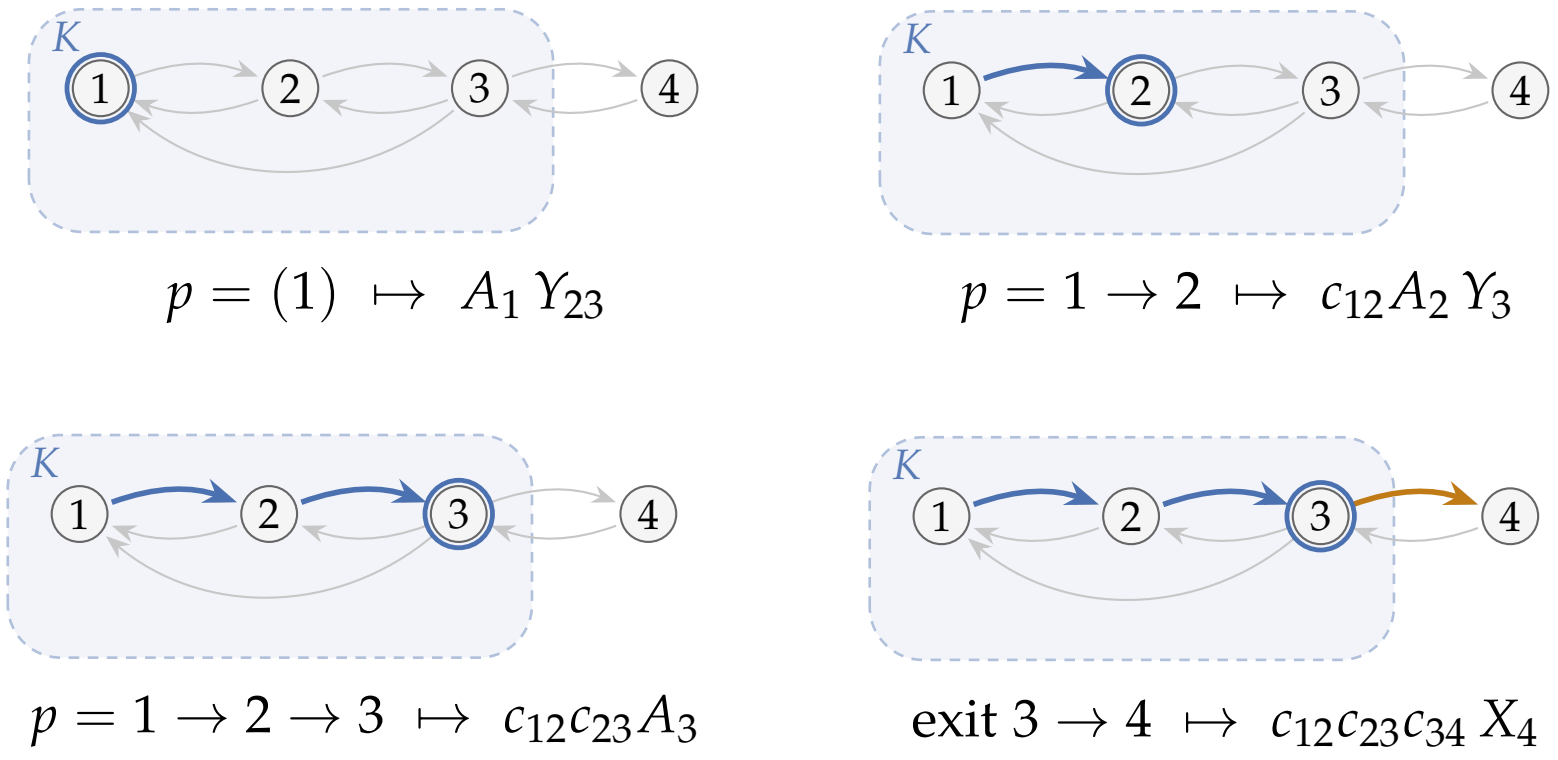
\captionof{figure}{The exit-path bijection for $K=\{1,2,3\}$, $i=1$. Following the arrows of an acyclic function on $K$ from $1$ traces a simple path $p$ that either stops at a root at $k$ (circled), leaving $A_kY_{K\setminus V(p)}$, or leaves $K$ along $k\to j$, leaving $X_j$. Conversely, a path, its final choice and an acyclic function on $K\setminus V(p)$ reassemble into an acyclic function on $K$.}\label{fig:exitpath}
\end{minipage}\end{center}
\medskip

\end{example}

Every term in \eqref{eq:exitpath} has at most one $X$-variable, indexed outside $K$, and a minor supported in $K\setminus\{i\}$. Thus $X_iY_K$ is already straightened. With more factors, we must keep any new $X$-variable outside the other minors as well.

\begin{theorem}[Weighted straightening]\label{thm:full}
Every monomial in the $X_i$ and weighted principal minors $Y_I$ is a linear combination of compatible monomials with coefficients in $\NN[c][A_1,\ldots,A_n]$.
\end{theorem}

\begin{proof}
First straighten the $Y$-factors. It remains to treat $X^aB_U$. Measure its incompatibility by
\[
d(a,U)=\sum_{i\in\supp U}a_i.
\]
If this is zero, there is nothing to do. Otherwise choose $i\in\supp U$ with $a_i>0$ and put $J=\supp U$. Apply the exit-path identity to $X_iY_J$, using $B_U=Y_JB_{U-\one_J}$.

All resulting minors are supported in $J$. Straighten them into $B_W$ with $\supp W\subseteq J$; Theorem~\ref{thm:B} never enlarges support. The exponent of $X_i$ has dropped by one, and any new $X_j$ lies outside $J$, hence outside $\supp W$. Thus $d$ decreases. Induction finishes the proof.
\end{proof}

This is why we use the outermost layer. Every remaining minor stays inside it, while any new $X$-variable lies outside it. Later straightening cannot bring that variable back into a $Y$-factor.

\begin{example}\label{ex:weighted-straightening}
Apply Example~\ref{ex:exitpath} to $X_1Y_{123}Y_{12}$. Then straighten $Y_{23}Y_{12}$ and use
$Y_3Y_{12}=Y_{123}+f_{23}Y_1+t$. The result is
\[
\begin{aligned}
X_1Y_{123}Y_{12}={}&
 A_1Y_{123}Y_2+tA_1Y_2+f_{12}f_{23}A_1\\
&+c_{12}A_2Y_{123}+c_{12}f_{23}A_2Y_1+c_{12}tA_2\\
&+c_{12}c_{23}A_3Y_{12}
+c_{12}c_{23}c_{34}X_4Y_{12}.
\end{aligned}
\]
All terms are compatible. In particular, the new $X_4$ is outside the remaining minor $Y_{12}$.
\end{example}

\subsection{The linear LP algebra}
Set $c_{ij}=1$. The weighted relations become the defining relations of $\cA_\Gamma$, and compatible monomials become cluster monomials \cite{LP2,DTW}. Theorem~\ref{thm:full} therefore proves positivity of multiplication in the cluster-monomial basis.

\begin{corollary}\label{cor:LP-positivity}
Every monomial in the cluster variables of $\cA_\Gamma$ is a linear combination of cluster monomials with coefficients in $\NN[A_1,\ldots,A_n]$. In particular, \cite[Conjecture~7.3(ii)]{LP2} holds for every directed graph.
\end{corollary}

\section{Laurent positivity}\label{sec:laurent}

We next fix a chart and expand every variable in its coordinates. Now we must choose our use of the pair rule so that we divide only by chart monomials. The nested structure provides the choice.

\subsection{The nested forest of a chart}
A \emph{nested collection} $\cN$ consists of nonempty strongly connected subsets of $[n]$ with the following properties:
\begin{enumerate}[leftmargin=2.2em,label=\textup{(\roman*)}]
\item Two members are nested or disjoint.
\item Any subcollection of pairwise disjoint members is exactly the set of strongly connected components of its union.
\end{enumerate}
Its support is $U=\bigcup_{I\in\cN}I$. It is \emph{maximal on $U$} if no further strongly connected subset of $U$ can be added. The corresponding chart is
\[
\mathbf x_{\cN}=\{Y_I:I\in\cN\}\cup\{X_j:j\notin U\}.
\]
At $c_{ij}=1$, these are precisely the clusters of $\cA_\Gamma$ \cite[Theorem~1.1]{LP2}.

In a maximal nested collection, each member $I$ has one \emph{private vertex}: the unique vertex of $I$ not contained in any smaller member. Equivalently, the map sending $v\in U$ to the smallest member containing $v$ is a bijection from $U$ to $\cN$ \cite[Lemma~3.2]{LP2}. The maximal members are the strongly connected components of $\Gamma[U]$. The \emph{children} of $I$ are its maximal proper members. These inclusions form a rooted forest, with each member consisting of its children and one private vertex. We will induct along this forest.

Let
\[
\cL_\cN=\NN[c,A_1,\ldots,A_n]
 [Y_I^{\pm1}:I\in\cN,\ X_j^{\pm1}:j\notin U].
\]
For now, this denotes the positive Laurent expressions in the chart elements, viewed as rational functions in our original variables. Algebraic independence of the chart will follow below. The set $\cL_\cN$ is closed under sums, products, and division by chart monomials.

\begin{theorem}[Weighted Laurent positivity]\label{thm:laurent}
Every $Y_S$ and every $X_i$ belongs to $\cL_\cN$.
\end{theorem}

\subsection{Removing a private vertex}
For a strongly connected set $D$, we first locate the smallest chart member containing it. Its private vertex will be the vertex we remove.

\begin{lemma}\label{lem:intrinsic}
Let $D\subseteq U$ be nonempty and strongly connected. There is a smallest member $I_D\in\cN$ containing $D$. If $m_D$ is its private vertex, then $m_D\in D$. Moreover,
\[
C_D=I_D\setminus\{m_D\}
\]
is the union of the children of $I_D$, and those children are the strongly connected components of $C_D$. In particular, $Y_{C_D}$ is a chart monomial.
\end{lemma}

\begin{proof}
The set $D$ lies in one maximal member, and the members containing it form a chain. The children of $I_D$ partition $I_D\setminus\{m_D\}$: every other vertex belongs to a proper member, while $m_D$ belongs to none. By nestedness, these children are the strongly connected components of their union. If $m_D\notin D$, strong connectivity would force $D$ into one child, contradicting the choice of $I_D$.
\end{proof}

The natural pair is now $(D,C_D)$: we can divide by the chart monomial $Y_{C_D}$, and a nonempty boundary heap must consume the private vertex.

\begin{proposition}[Private-vertex identity]\label{prop:private}
Suppose $D\subseteq U$ is strongly connected and $D\notin\cN$. Put $I=I_D$, $m=m_D$, $C=C_D$, and $\kappa_D=\one_D+\one_C$. Then
\begin{equation}\label{eq:private}
Y_DY_C=Y_IY_{D\setminus\{m\}}
+\sum_{\substack{H\ne\emptyset\\H\text{ a $(D,C)$-boundary heap}}}
 c^H B_{\kappa_D-\cont(H)}.
\end{equation}
Every heap in the sum uses $m$ exactly once. Both residual layers are therefore contained in $C$.
\end{proposition}

\begin{proof}
Apply the pair rule to $(D,C)$. Since $m\in D$, their union is $I$ and their intersection is $D\setminus\{m\}$, giving the empty-heap term. Every minimal piece of a nonempty boundary heap must leave $C$, so it passes through $m$. But $m$ has capacity one. It is used exactly once and does not occur in the residual vector.
\end{proof}

\begin{example}\label{ex:private}
Take the chain chart
\[
\cN=\{1,12,123,1234\},
\qquad\mathbf x_\cN=(Y_1,Y_{12},Y_{123},Y_{1234}).
\]
For $D=23$, we have $I_D=123$, $m_D=3$, and $C_D=12$. The private-vertex identity is
\[
Y_{23}Y_{12}=Y_{123}Y_2+tY_2+f_{12}f_{23}.
\]
The heaps for the ordered pair $(23,12)$ are $\emptyset$, $\tau$, and $\beta\alpha$. The nonempty ones consume vertex $3$. After division by $Y_{12}$, only the smaller minor $Y_2$ remains to be expanded.

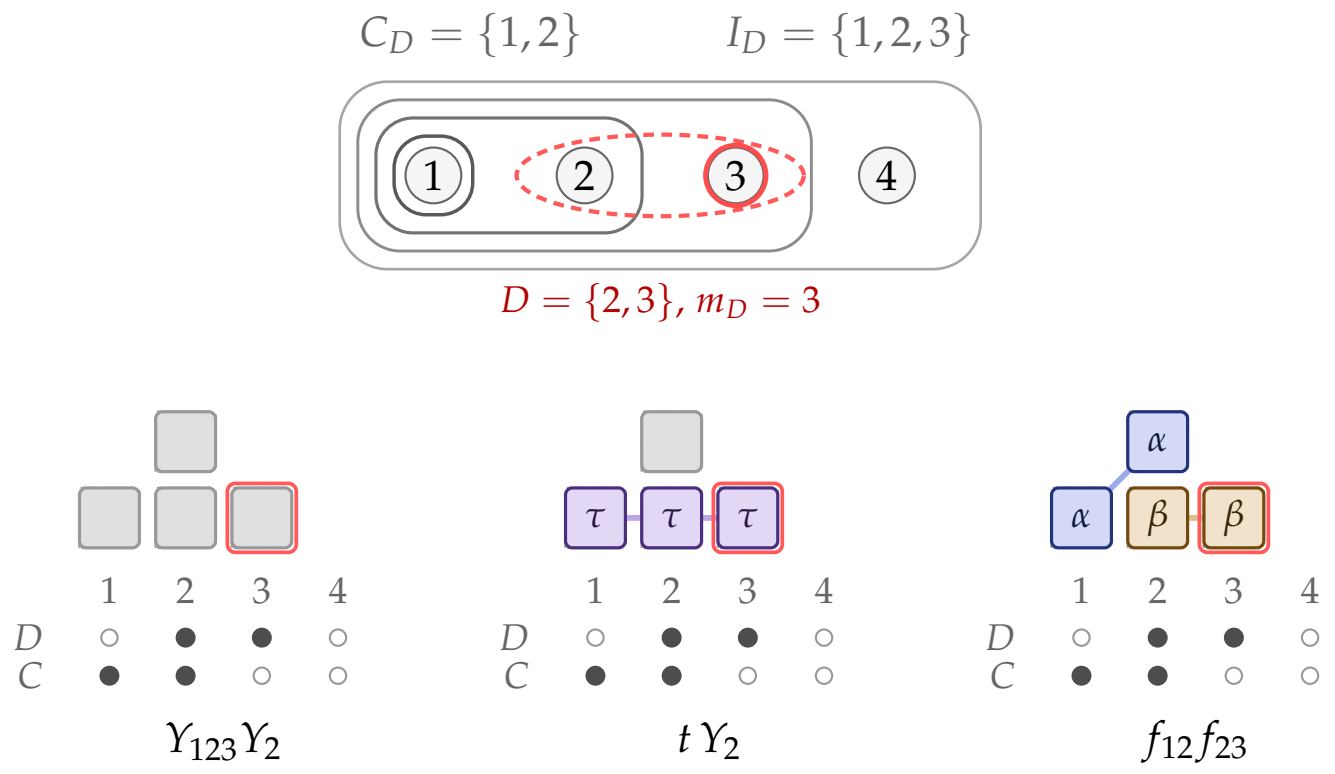
\begin{figure}[!tbp]
\centering
\begin{tikzpicture}[baseline=(base)]
  \coordinate (base) at (0,-0.1);
  \draw[black!35,rounded corners=9pt,line width=.5pt] (-0.62,-0.62) rectangle (3.62,0.62);
  \draw[black!45,rounded corners=8pt,line width=.55pt] (-0.5,-0.5) rectangle (2.5,0.5);
  \draw[black!55,rounded corners=7pt,line width=.6pt] (-0.38,-0.38) rectangle (1.38,0.38);
  \draw[black!65,rounded corners=6pt,line width=.65pt] (-0.26,-0.26) rectangle (0.26,0.26);
  \draw[red!65,dash pattern=on 2pt off 1.4pt,line width=.8pt] (1.5,0) ellipse (0.95 and 0.27);
  \node[gvertsm] (v1) at (0,0) {$1$};
  \node[gvertsm] (v2) at (1,0) {$2$};
  \node[gvertsm] (v3) at (2,0) {$3$};
  \node[gvertsm] (v4) at (3,0) {$4$};
  \draw[red!70,line width=.9pt] (v3.center) circle (0.2);
  \node[font=\scriptsize,text=black!60] at (2.75,0.93) {$I_D=\{1,2,3\}$};
  \node[font=\scriptsize,text=black!60] at (0.25,0.93) {$C_D=\{1,2\}$};
  \node[font=\tiny,text=red!70!black] at (1.5,-0.84) {$D=\{2,3\}$,\ $m_D=3$};
\end{tikzpicture}
\\[12pt]
\newcommand\privpic[2]{\begin{tikzpicture}[baseline=0pt,scale=1.2]
  \cslot{1}{1}\cslot{2}{1}\cslot{2}{2}\cslot{3}{1}#1
  \draw[red!65,line width=.7pt,rounded corners=1.5pt] (0.86,0.02) rectangle (1.24,0.4);
  \ccols{4}
  \membrow{0.42}{4}{-0.45}{$D$}{2,3}\membrow{0.42}{4}{-0.66}{$C$}{1,2}
  \node[font=\scriptsize] at (0.84,-1.02) {#2};
\end{tikzpicture}}%
\privpic{\ccellf{1}{1}{black!12}{black!40}{}\ccellf{2}{1}{black!12}{black!40}{}\ccellf{2}{2}{black!12}{black!40}{}\ccellf{3}{1}{black!12}{black!40}{}}{$Y_{123}Y_2$}\hspace{18pt}
\privpic{\cjoin{1}{1}{3}{1}{ctau}\ccell{1}{1}{ctau}{\tau}\ccell{2}{1}{ctau}{\tau}\ccell{3}{1}{ctau}{\tau}\ccellf{2}{2}{black!12}{black!40}{}}{$t\,Y_2$}\hspace{18pt}
\privpic{\cjoin{2}{1}{3}{1}{cbeta}\cjoin{1}{1}{2}{2}{calpha}\ccell{2}{1}{cbeta}{\beta}\ccell{3}{1}{cbeta}{\beta}\ccell{1}{1}{calpha}{\alpha}\ccell{2}{2}{calpha}{\alpha}}{$f_{12}f_{23}$}
\captionof{figure}{The private-vertex identity for $D=\{2,3\}$ in the chart $\{1\}\subset\{1,2\}\subset\{1,2,3\}\subset[4]$. Top: the chart as nested boxes, with $D$ dashed and its private vertex $m_D=3$ circled. Below: the three $(D,C_D)$-boundary heaps in the capacity $\kappa_D=\one_D+\one_{C_D}$. Every nonempty heap fills the single box over $m_D$ (framed), so the leftover layers lie inside the child $C_D=\{1,2\}$.}\label{fig:private}
\end{figure}

\end{example}

\subsection{The three-step proof}
\begin{proof}[Proof of Theorem~\ref{thm:laurent}]
\emph{Minors supported in $U$.}
Factor over strongly connected components. For a strongly connected $D$, induct on $|I_D|$. If $D\in\cN$, its minor is a chart variable. Otherwise divide \eqref{eq:private} by the chart monomial $Y_{C_D}$. Apart from $Y_{I_D}$, every factor is supported in $C_D$. Its strongly connected components lie in children of $I_D$, so induction applies. All operations are subtraction-free.

\emph{The variables $X_i$.}
If $i\notin U$, the variable is already in the chart. If $i\in U$, apply the exit-path identity with $K=U$ and divide by $Y_U$, a chart monomial. The minors on the right were handled in the first step, and any $X_j$ on the right lies outside $U$.

\emph{Arbitrary minors.}
Induct on $|S\setminus U|$. If this is positive, choose $i\in S\setminus U$, apply the exit-path identity to $X_iY_S$, and divide by the chart variable $X_i$. Every path removes $i$, so the remaining minors have fewer vertices outside $U$. The $X$-variables were handled in the second step.
\end{proof}

\begin{corollary}\label{cor:chart}
The $n$ elements of a chart are algebraically independent over
$\QQ(c)(A_1,\ldots,A_n)$. Their Laurent expansions are unique, and changes between charts are subtraction-free birational maps.
\end{corollary}

\begin{proof}
The field $F$ has transcendence degree $n$ over the stated coefficient field. Theorem~\ref{thm:laurent} expresses every $X_i$ rationally in the chart, so its $n$ elements generate $F$. They must therefore be algebraically independent.
\end{proof}

\begin{corollary}[Laurent positivity for linear LP algebras]\label{cor:laurent-LP}
At $c_{ij}=1$, every cluster variable is a Laurent polynomial with coefficients in $\NN[A_1,\ldots,A_n]$ in every cluster. Every exchange polynomial also has nonnegative coefficients. Thus \cite[Conjecture~7.1]{LP2} holds; Laurent positivity of the cluster variables was proved in \cite{Glew}.
\end{corollary}

\begin{proof}
The charts become clusters, so the first assertion follows. If mutation replaces $x$ by $x'$, then the exchange Laurent polynomial is $xx'$, which is positive in the original cluster. Its polynomial normalization differs by a monomial in the other cluster variables, so its coefficients remain nonnegative.
\end{proof}

\begin{remark}[Relation to the first author's proof]\label{rem:glew}
The three steps above follow the strategy of \cite[Theorem~5.1]{Glew}. There, a positive path identity, obtained by iterating determinantal identities, drives an induction along an activation order. Here Theorem~\ref{thm:B} gives the positive correction in \eqref{eq:private}, and the first induction follows the nested forest itself. Both corrections have the same sum, but their positive summands need not agree term by term. The intrinsic recursion avoids any choice of activation order.
\end{remark}

\begin{remark}[Iterating the recursion]\label{rem:laurent-expansions}
The proof of Theorem~\ref{thm:laurent} also tells us how to compute the expansions. For a minor supported in $U$, repeatedly apply the private-vertex identity and expand the smaller minors in the residual layers. Each resulting term records the boundary heaps chosen along the way. Exit paths handle the remaining variables. All denominators are chart monomials, so no subtraction is needed at any stage.
\end{remark}

\begin{example}[A positive Laurent expansion]\label{ex:laurent-iteration}
Return to the chain chart
$\cN=\{1,12,123,1234\}$ of Example~\ref{ex:private}.
The first private-vertex identity gives
\[
Y_{12}Y_{23}=Y_{123}Y_2+tY_2+f_{12}f_{23}.
\]
To expand the remaining minor $Y_2$, apply the identity once more,
now at the private vertex $2$ of $12$:
\[
Y_1Y_2=Y_{12}+f_{12}.
\]
Substituting and clearing denominators gives
\begin{equation}\label{eq:Y23}
\boxed{
Y_{23}=\frac{Y_{12}Y_{123}+f_{12}Y_{123}
+tY_{12}+tf_{12}+f_{12}f_{23}Y_1}{Y_1Y_{12}}.
}
\end{equation}
There are three choices at the first step. Two leave a factor $Y_2$, which gives two terms at the next step; the third leaves no minor to expand. This accounts for the five positive terms. Nothing special about this example is needed for the procedure: the private vertex always reduces the problem to smaller chart members, and exit paths handle variables outside the support.
\end{example}


\begin{thebibliography}{BCKZ2}
\bibitem[BCKZ1]{BCKZ1} E.~Banaian, S.~Chepuri, E.~Kelley, S.~W.~Zhang, \emph{Rooted clusters for graph LP algebras}, SIGMA 18 (2022), 089; arXiv:2107.14785.
\bibitem[BCKZ2]{BCKZ2} E.~Banaian, S.~Chepuri, E.~Kelley, S.~W.~Zhang, \emph{Snake graphs for graph LP algebras}, arXiv:2312.12313.
\bibitem[CF]{CF} P.~Cartier, D.~Foata, \emph{Probl\`emes combinatoires de commutation et r\'earrangements}, Lecture Notes in Math.\ 85, Springer, 1969.
\bibitem[Ch]{Chaiken} S.~Chaiken, \emph{A combinatorial proof of the all minors matrix tree theorem}, SIAM J.\ Algebraic Discrete Methods 3 (1982), 319--329.
\bibitem[DTW]{DTW} G.~Z.~Dantas e Moura, R.~C.~Telekicherla Kandalam, D.~Woodruff, \emph{Cluster monomials in graph Laurent phenomenon algebras}, Algebraic Combinatorics 8 (2025), 997--1019.
\bibitem[FZ1]{FZ1}
S.~Fomin and A.~Zelevinsky,
\emph{Cluster algebras I: Foundations},
J.~Amer.~Math.~Soc. 15 (2002), 497--529.
\bibitem[GM]{GarsiaMilne} A.~M.~Garsia, S.~C.~Milne, \emph{Method for constructing bijections for classical partition identities}, Proc.\ Natl.\ Acad.\ Sci.\ U.S.A.\ 78 (1981), 2026--2028.
\bibitem[G]{Glew} R.~Glew, \emph{Positivity for graph Laurent phenomenon algebras}, arXiv:2610.08914.
\bibitem[GHKK]{GHKK}
M.~Gross, P.~Hacking, S.~Keel, and M.~Kontsevich,
\emph{Canonical bases for cluster algebras},
J.~Amer.~Math.~Soc. 31 (2018), 497--608.
\bibitem[LP1]{LP1} T.~Lam, P.~Pylyavskyy, \emph{Laurent phenomenon algebras}, Cambridge J.\ Math.\ 4 (2016), 121--162; arXiv:1206.2611.
\bibitem[LP2]{LP2} T.~Lam, P.~Pylyavskyy, \emph{Linear Laurent phenomenon algebras}, Int.\ Math.\ Res.\ Not.\ IMRN (2016), 3163--3203; arXiv:1206.2612.
\bibitem[LS]{LeeSchiffler}
K.~Lee and R.~Schiffler,
\emph{Positivity for cluster algebras},
Ann.~of~Math. 182 (2015), 73--125.
\bibitem[LyP]{LyonsPeres} R.~Lyons, Y.~Peres, \emph{Probability on Trees and Networks}, Cambridge University Press, 2016.
\bibitem[PW]{PW} J.~G.~Propp, D.~B.~Wilson, \emph{How to get a perfectly random sample from a generic Markov chain and generate a random spanning tree of a directed graph}, J.\ Algorithms 27 (1998), 170--217.
\bibitem[T]{Tang}
Q.~Tang, \emph{Positive bases and polytopes for Laurent
phenomenon algebras},
Hexagon preprint hexagon:2611.00284v1 (2026);
\url{https://hexagonmath.org/2611.00284v1}.
\bibitem[V]{Viennot} G.~X.~Viennot, \emph{Heaps of pieces, I: Basic definitions and combinatorial lemmas}, in: Combinatoire \'enum\'erative, Lecture Notes in Math.\ 1234, Springer, 1986, 321--350.
\bibitem[Wi]{Williams}
L.~K.~Williams,
\emph{Cluster algebras: An introduction},
Bull.~Amer.~Math.~Soc. 51 (2014), 1--26.
\bibitem[W]{Wilson} D.~B.~Wilson, \emph{Generating random spanning trees more quickly than the cover time}, in: Proceedings of the 28th Annual ACM Symposium on Theory of Computing, 1996, 296--303.
\end{thebibliography}
\end{document}